\documentclass[11pt, a4paper]{article}

\usepackage{xcolor}

\definecolor{rose}{RGB}{207, 96, 122}

\definecolor{kit-green}{RGB}{0, 150, 130}
\colorlet{kit-green100}{kit-green}
\colorlet{kit-green90}{kit-green!90!white}
\colorlet{kit-green80}{kit-green!80!white}
\colorlet{kit-green70}{kit-green!70!white}
\colorlet{kit-green60}{kit-green!60!white}
\colorlet{kit-green50}{kit-green!50!white}
\colorlet{kit-green40}{kit-green!40!white}
\colorlet{kit-green30}{kit-green!30!white}
\colorlet{kit-green25}{kit-green!25!white}
\colorlet{kit-green20}{kit-green!20!white}
\colorlet{kit-green15}{kit-green!15!white}
\colorlet{kit-green10}{kit-green!10!white}
\colorlet{kit-green5}{kit-green!5!white}

\definecolor{kit-royalblue}{RGB}{0, 45, 76}
\colorlet{kit-royalblue100}{kit-royalblue}
\colorlet{kit-royalblue90}{kit-royalblue!90!white}
\colorlet{kit-royalblue80}{kit-royalblue!80!white}
\colorlet{kit-royalblue70}{kit-royalblue!70!white}
\colorlet{kit-royalblue60}{kit-royalblue!60!white}
\colorlet{kit-royalblue50}{kit-royalblue!50!white}
\colorlet{kit-royalblue40}{kit-royalblue!40!white}
\colorlet{kit-royalblue30}{kit-royalblue!30!white}
\colorlet{kit-royalblue25}{kit-royalblue!25!white}
\colorlet{kit-royalblue20}{kit-royalblue!20!white}
\colorlet{kit-royalblue15}{kit-royalblue!15!white}
\colorlet{kit-royalblue10}{kit-royalblue!10!white}
\colorlet{kit-royalblue5}{kit-royalblue!5!white}

\usepackage[utf8]{inputenc}
\usepackage[T1]{fontenc}
\usepackage{lmodern}
\usepackage{geometry}
\usepackage{setspace}
\usepackage{parskip}
\usepackage{fancyhdr}
\usepackage{titlesec}
\usepackage{graphicx}
\usepackage{hyperref}
\usepackage{abstract}
\usepackage{blindtext}  
\usepackage{amsmath}
\usepackage{amssymb}
\usepackage{bm}
\usepackage{mathrsfs}
\usepackage{subcaption}
\usepackage[svgnames]{xcolor}
\usepackage{float}
\usepackage{adjustbox}
\usepackage{csquotes}
\usepackage{tikz}
\usepackage{tabularx}
\usepackage{pgfplots}
\usepackage[normalem]{ulem}
\usepackage{booktabs}
\pgfplotsset{compat=1.18}

\hypersetup{
    colorlinks=true,
    linkcolor=rose,
    urlcolor=black,
    citecolor=rose,
}

\titleformat{\section}{\large\bfseries}{\thesection}{1em}{}
\titleformat{\subsection}{\normalsize\bfseries}{\thesubsection}{1em}{}

\renewcommand{\vec}[1]{\bm{#1}}  
\renewcommand{\d}{\mathrm{d}}
\newcommand{\npe}{^{n+1}}  
\newcommand{\npeh}{^{n+\frac{1}{2}}}

\newcommand{\n}{^{n}}  
\newcommand{\transp}{^\top}  

\renewcommand{\epsilon}{\varepsilon}

\usepackage{amsthm}
\usepackage{cleveref}

\theoremstyle{plain}
\newtheorem{theorem}{Theorem}[section]
\crefname{theorem}{theorem}{theorems}
\Crefname{theorem}{Theorem}{Theorems}

\theoremstyle{definition}

\crefname{definition}{definition}{definitions}
\Crefname{definition}{Definition}{Definitions}

\theoremstyle{remark}
\newtheorem{remark}[theorem]{Remark}
\crefname{remark}{remark}{remarks}
\Crefname{remark}{Remark}{Remarks}

\newenvironment{authcontrib}[1]{%
    \subsection*{\textnormal{\textbf{Author Contributions}}}%
    \noindent #1}%
{}%
\newenvironment{acks}[1]{%
    \subsection*{\textnormal{\textbf{Acknowledgements}}}%
    \noindent #1}%
{}%
\newenvironment{dci}[1]{%
    \subsection*{\textnormal{\textbf{Declaration of competing interests}}}%
    \noindent #1}%
{}%
\newenvironment{code}[1]{%
    \subsection*{\textnormal{\textbf{Code}}}%
    \noindent #1}%
{}%

\usepackage{dashbox}%

\newcommand\colordashedph[1][H]{\setlength{\fboxsep}{0.3pt}\setlength{\dashlength}{2.2pt}\setlength{\dashdash}{1.1pt} \dbox{\colorbox{lightgray!40!Lavender!40!}{\phantom{#1}}}}

\renewcommand{\square}{\colordashedph}
\usepackage{ifthen}
\usepackage[]{changes}
\setaddedmarkup{\color{blue}\uuline{#1}}
\setdeletedmarkup{\color{red}\sout{#1}}
\definechangesauthor[color=darkgray]{pb}
\definechangesauthor[color=darkgray]{plk}
\definechangesauthor[color=darkgray]{ll}

\newif\ifshowchanges
\showchangestrue    

\usepackage{marginnote}

\usepackage[backend=biber, style=numeric-comp, giveninits=true, maxbibnames=4, sortcites=true, sorting=none, url=true, isbn=false, date=year]{biblatex}
\title{Port-Hamiltonian multibody dynamics:
Lagrangian formulation, consistent interconnection, structure-preserving simulation and index-reduction
}

\author{
\begin{tabular}{ccccc}
Lisa Latussek\footnote{Karlsruhe Institute of Technology (KIT)} \textsuperscript{,}\footnote{ETH Zurich} & &
Philipp L. Kinon\footnotemark[1] \textsuperscript{,}\footnote{Eindhoven University of Technology (TU/e)} & &  Peter Betsch\footnotemark[1] \\
\href{mailto:llatussek@ethz.ch}{\texttt{\small llatussek@ethz.ch}} &\quad& \href{mailto:philipp.kinon@kit.edu}{\texttt{\small philipp.kinon@kit.edu}} &\quad& \href{mailto:peter.betsch@kit.edu}{\texttt{\small peter.betsch@kit.edu}}
\end{tabular}
}

\date{\today}

\begin{document}

\maketitle

\begin{abstract}
    \noindent
    This work introduces a port-Hamiltonian (PH) model for constrained mechanical systems, which is directly derived from the Lagrangian equations of motion. The present PH framework incorporates a singularity-free director representation of rigid body rotations, resulting in constant mass matrices.
It is shown that the power-preserving interconnection of PH rigid-body subsystems is mathematically equivalent to the classical description of ideal joints using kinematic pairs. This establishes a PH multibody dynamics framework that is consistent with traditional modeling paradigms.

Notably, the PH structure of the governing index-2 differential-algebraic equations enables the application of an implicit, structure preserving midpoint time integration. 
The proposed scheme is able to satisfy both the balance laws for total energy and angular momentum as well as the position-level constraints.
These properties make the proposed method remarkably robust and enable
stable long-term simulations. 
Furthermore, a variationally derived index-reduction strategy is incorporated that enforces velocity-level constraints in addition to position-level constraints while preserving the port-Hamiltonian structure. 
Numerical examples illustrate the favorable properties of the proposed formulation, which is well-suited for energy-based control design.

\vspace{0.5em}
{\bfseries Keywords:}\hspace{0.8em} Multibody dynamics \(\cdot\) Director formulation \(\cdot\) Port-Hamiltonian systems \(\cdot\) Energy-based methods \(\cdot\) Structure-preserving integration \(\cdot\) Index reduction
\par
\end{abstract}


\section{Introduction}\label{sec:introduction}

The dynamic behavior of multibody systems (MBS), described as sets of rigid bodies interconnected by joints
restricting relative motion \cite{simeon_2013_computational}, can become highly intricate, involving e.g. large rotations or interactions across multiple physical domains. Accurate and robust simulation of
such systems is essential, particularly for long-term studies, where
structure-preserving numerical schemes typically exhibit superior behavior \cite{hairer2006}.

A promising approach to addressing these challenges is the use of the \emph{port-Hamiltonian} (PH) framework, which has emerged as an advantageous modeling paradigm for complex dynamical systems across
multiple research fields \cite{duindam2009,vanderschaft2014}.
A key benefit of the PH representation lies in its explicit formulation of \emph{ports} that provide
interfaces for energy exchange between subsystems. This allows modular interconnection of
submodules while ensuring that global power balances are preserved. This property has also
proven useful in the context of (flexible) multibody dynamics \cite{macchelli2009,lohmayer_2024_exergetic,berger2025,kinonthoma2023,brugnoli2021},
as it enables
energy-consistent
coupling of individual bodies. In addition, external ports provide an ideal
interface for control design , see e.g. \cite{ortega_2002_interconnection,caasenbrood_2022_energyshaping,vanderschaft2000L2}.

Remarkably, many previous works formulated in the PH framework are
based on positions and momenta \cite{vanderschaft2014,macchelli2009,duindam_2009_modeling,forni2015,trivedi_2011_modeling,ayala_2023_energybased,warsewa_2021_porthamiltonian} and not positions and velocities (corresponding to the Lagrangian formulation of mechanics). However, the Lagrangian description
is oftentimes the preferred choice
\cite{curiel_2014_classical,leyendecker_2008_variational,bauchau_2009_scaling,gerstmayr_2006_3d,bayo_1996_augmented,laulusa_2008_review,jain_1995_diagonalized,gay-balmaz_2018_lagrangian,chung_2009_cooperative,betsch2001, betsch2002, betsch2006}.

While the treatment of rigid body rotations in the context of the PH framework has so far been approached
using Lie-group formulations \cite{macchelli2009,lohmayer_2024_exergetic}, coordinate-based approaches can also be a suitable option. A first step in this direction
has been made in \cite{forni2015}. Formulations of MBS relying on minimal coordinates introduce singularities and can result in numerical instabilities,
particularly in the presence of large rotations. Alternatives are given by redundant coordinates for the description of finite rotations
such as
unit-quaternions \cite{moller_2012_rigid,harsch_2023_nonunit,kinon_2024_conserving,may_2025_galerkinbased} or directors \cite{krenkConservativeRigidBody2014,bauchau_2011_flexible,betsch_2003_constrained}, see overview works \cite{bauchau_2014_interpolation,romero_2004_interpolation,ibrahimbegovic_1997_choice}.
The latter approach proves to be particularly attractive, since, unlike quaternion representations,
constant mass matrices are obtained. 
Furthermore, the respective internal constraints of rigidity are only quadratic \cite{betsch2001,betsch2002}. However, this comes at the price of six additional Lagrange multipliers enforcing the above-mentioned constraints.

In the director formulation of rigid body dynamics
, each body is equipped with a set of director vectors
that directly parameterize its orientation, avoiding the singularities inherent in
minimal-coordinate representations.
The director formulation
has been extended to include external wrenches \cite{betsch2013} and kinematic pairs \cite{betsch2006,sanger2011}.
This enables a systematic composition of complex multibody systems.

The motion of constrained mechanical systems is governed by \emph{differential-algebraic equations} (DAEs) with differential index three.
A PH formulation that naturally includes DAE systems has been proposed in \cite{mehrmann19}.
From a computational point of view, it is 
desirable to obtain DAE representations 
with a lower index \cite{simeon_2013_computational,gear_1988_differentialalgebraic,mattsson_1993_index,kunkel_2004_index}, as they are less prone to numerical instabilities. A well-known approach here is the GGL stabilization \cite{gear_1985_automatic}, where additionally velocity level constraints are considered. While the original GGL stabilization 
sacrifices the Hamiltonian structure of MBS, thus hindering an appealing structure-preserving discretization, a variational justification and modification called the \textit{GGL principle} has been recently developed \cite{kinon2023structure-preserving,kinon2023the-ggl-variational}.

\subsection{Research gap}\label{sec:research_gap}
Previous works have addressed PH formulations for multibody systems \cite{macchelli2009,brugnoli2021,berger2025} and 
rigid body dynamics 
\cite{forni2015}, but a systematic combination of the director
formulation \cite{betsch2001}, PH framework, and structure-preserving discretization for 
MBS is still missing. This work pursues that connection. Embedding the director formulation in the PH framework should ensure energy consistency, modular interconnection, and a natural setting for control; and structure-preserving time-stepping
enables the conservation of 
invariants such as energy and angular momentum over long horizons in discrete time \cite{hairer2006}. Further, enhanced numerical stability 
can be gained through index-reduction using the GGL principle \cite{kinon2023the-ggl-variational}.

\subsection{Contributions}\label{sec:contributions}
The contributions of this work can be summarized as follows:
\begin{itemize}
    \item[\textbf{C1}] We present a new PH-DAE formulation for mechanical systems based on the Lagrangian equations of the first kind and specify this formulation for multibody systems.
    \item[\textbf{C2}] We introduce a novel PH-DAE formulation for multibody systems based on the director formulation, which accommodates large rotations and guarantees intrinsic energy conservation. This formulation enables seamless
          interconnection with other PH systems and benefits from a constant mass matrix.
    \item[\textbf{C3}] We analyze the structure-preserving time discretization of the proposed framework using the implicit midpoint rule. We show that this approach preserves the balance of energy and angular momentum of
          the system, ensuring long-term accuracy and stability. The scheme also prevents drift-off from position-level constraints.
    \item[\textbf{C4}] We demonstrate that the port matrices describing the interconnection of rigid bodies as PH systems encode the same information as kinematic pairs in classical multibody dynamics, providing a rigorous link between the
          two perspectives.
    \item[\textbf{C5}] We achieve an index-reduction within the PH-DAE framework using the GGL principle leading to an additional satisfaction of velocity-level constraints in discrete time.
\end{itemize}

\subsection{Notation}\label{annex_notation}
We follow the following notation rules.
Scalars are represented in italic fonts, e.g. $f,g \in \mathbb{R}$.
Vectors from the three-dimensional Euclidean space, e.g. $\vec{a},\vec{b},\vec{c} \in \mathbb{E}^3$, and tuples $\vec{e} \in \mathbb{R}^3$
are treated equally by using bold letters. Thus,
each vector can be identified with its component tuple in the basis $\{\vec{e}_i\}_{i=1}^3$ of a right-handed, orthonormal, inertial frame
such that
\begin{equation*}
    \vec{x} = \sum_{i=1}^3 x_i \vec{e}_i = x_i \vec{e}_i \, \hat{=} \begin{bmatrix} x_1 \\ x_2 \\ x_3 \end{bmatrix},
\end{equation*}
unless noted otherwise.
Note that the Einstein summation convention for repeated indices is used throughout this work, e.g. for $i,j,k \in \{1,2,3\}$.
A short bracket notation is used to express $\vec{x}=[\vec{a}\transp,\vec{b}\transp]\transp =:(\vec{a},\vec{b})$.
The scalar product of two vectors is denoted by e.g. $\vec{a}\transp \vec{b} = a_i b_i$ and the cross-product of two vectors $\vec{c} = \vec{a} \times \vec{b} = \epsilon_{ijk} a_j b_k \vec{e}_i = \hat{\vec{a}}\vec{b} = -\hat{\vec{b}}\vec{a}$, where we made use of the well-known Levi-Civita symbol $\epsilon_{ijk}$ and introduced
\begin{align} \label{eq_skew_map}
    \hat{\square} = \begin{bmatrix}
                        0          & -\square_3 & \square_2  \\
                        \square_3  & 0          & -\square_1 \\
                        -\square_2 & \square_1  & 0
                    \end{bmatrix} \in \mathfrak{so}(3) , \quad \forall \ \square \in \mathbb{R}^3 ,
\end{align}
the mapping between $\mathbb{R}^3$ and the space of skew-symmetric matrices $\mathfrak{so}(3) \subset \mathbb{R}^{3 \times 3}$.
A matrix $\vec{A} \in \mathbb{R}^{n \times m}$ is denoted using bold face and its transpose reads
$\vec{A}\transp$.
Matrix-vector multiplication is denoted as $\vec{A} \vec{b}$.
Identity matrices are represented as $\vec{I}$ and matrices full of zeros as $\vec{0}$, their dimension being clear from the context.
The square $\square$ represents a placeholder.
The partial derivative of a function $f$ with respect to $\square$ is denoted $\partial_{\square} f$.
We introduce the special notation $\dot{\square} := \partial_t \square$ for time derivatives.
In this work gradients $\nabla g(\vec{q})$ and partial derivatives
\begin{align}
    \partial_{\vec{q}} f(\vec{q}, \vec{v})
    = \begin{bmatrix}
          \partial_{q_1}f(\vec{q}, \vec{v}) \\
          \partial_{q_2}f(\vec{q}, \vec{v}) \\
          \vdots                            \\
          \partial_{q_n}f(\vec{q}, \vec{v})
      \end{bmatrix}
\end{align}
are represented as column vectors.
Additionally, second partial derivatives yield matrices in such a way that
\begin{equation}
    \begin{aligned}
        \partial^2_{\vec{q} \vec{v}} f(\vec{q}, \vec{v}) :=
        \frac{\partial}{\partial \vec{v}}\left( \frac{\partial f(\vec{q}, \vec{v})}{\partial \vec{q}}  \right)
        = (\partial^2_{q_i v_j} f) \vec{e}_i \vec{e}_j\transp .
    \end{aligned}
\end{equation}
Due to Schwarz' theorem $\partial^2_{q_i v_j} f = \partial^2_{v_j q_i} f$, the identity $\partial^2_{\vec{q}\vec{v}} f(\vec{q}, \vec{v}) = \partial^2_{\vec{v}\vec{q}} f(\vec{q}, \vec{v})\transp$
holds true.

\subsection{Outline}\label{sec:outline}

The remainder of this work is structured as follows.
In Section~\ref{sec:PHS_primer} we briefly introduce the employed PH model.
Next, the general Lagrangian framework for mechanics is brought into the PH form in Section~\ref{sec_PH_lagrangian} and we discuss relevant special cases.
Section \ref{sec:modeling} introduces a multibody framework based on the director formulation for rigid bodies and kinematic pairs. Specifically, we establish
the link between kinematic pairs and PH power-preserving interconnections in a rigorous way. Section \ref{sec:numerical_integration} translates the framework to the discrete setting by applying the implicit midpoint rule and
demonstrating the exact representation of energy and angular momentum balances at the discrete level. In Section \ref{sec:GGL_principle}, we extend the framework by using a variational version of the well-known \textit{GGL stabilization}, the GGL principle, for the sake of index-reduction and the exact enforcement of constraints also on velocity-level.
Section \ref{sec:numerical_results} presents numerical
experiments on representative multibody systems, illustrating robustness, stability, and structure-preserving properties of the proposed framework. Finally, Section \ref{sec:conclusion} provides an outlook for future research, and summarizes the main
conclusions. Detailed derivations and supplementary material are included in the Appendices.

\section{Port-Hamiltonian systems}\label{sec:PHS_primer}
Let time $t \in [0, t_f]$ for some final time $t_f >0$. A \emph{port-Hamiltonian descriptor system} (PH-DAE) is a dynamical system characterized by a set of differential-algebraic equations
\cite{mehrmann19,kinon2025discrete,mehrmann_2023_control,beattie_2018_linear}
\begin{subequations} \label{eq:port_hamiltonian_form}
    \begin{align}
        \vec{E}(\vec{x}) \dot{\vec{x}} & = \vec{J}(\vec{x})\, \vec{z}(\vec{x}) + \vec{B}(\vec{x})\,  \vec{u}, \label{eq:state_evolution} \\
        \vec{y}                        & = \vec{B}(\vec{x})\transp \vec{z} . \label{eq:relate_costate_output}
    \end{align}
    together with a Hamiltonian function $H \in \mathcal{C}^1(\mathbb{R}^d,\mathbb{R})$,
    such that the
    constitutive relation linking the \emph{costate} function $ \vec{z}: \mathbb{R}^d \rightarrow \mathbb{R}^d$ and the Hamiltonian gradient,
    \begin{align}
        \vec{E}(\vec{x})\transp \vec{z}(\vec{x}) & = \nabla H(\vec{x}), \label{eq:relate_costate_Hamiltonian_differential_gradient}
    \end{align}
\end{subequations}
is verified for all \emph{state} vectors $\vec{x} \in \mathbb{R}^d $.

While the evolution of the state is governed by the differential(-algebraic) equations \eqref{eq:state_evolution}, relation \eqref{eq:relate_costate_output} captures the output behavior of the system. To this end, in \eqref{eq:port_hamiltonian_form}, $ \vec{u} \in \mathbb{R}^p $ represents \emph{input} quantities and $ \vec{y} \in \mathbb{R}^p $ denotes the power-conjugated and
collocated \emph{outputs}.
In addition, several matrices appear in \eqref{eq:port_hamiltonian_form}. These are the symmetric and state-dependent \emph{descriptor matrix} $\vec{E}(\vec{x})=\vec{E}(\vec{x})\transp \in \mathbb{R}^{d \times d}$, the skew-symmetric and state-dependent \emph{structure matrix} $\vec{J}(\vec{x})=-\vec{J}(\vec{x})\transp \in \mathbb{R}^{d \times d}$, and the state-dependent \emph{port matrix} $\vec{B}(\vec{x}) \in \mathbb{R}^{d \times p}$.
In the presence of algebraic constraints in \eqref{eq:state_evolution}, the descriptor matrix $\vec{E}$ is explicitly allowed to be singular.

By design then, PH-DAE systems ensure energetic consistency due to
the skew-symmetry of $\vec{J}(\vec{x})$ and relation \eqref{eq:relate_costate_Hamiltonian_differential_gradient}. For this, we consider the total time derivative of the Hamiltonian.
Making use of \eqref{eq:state_evolution} and \eqref{eq:relate_costate_Hamiltonian_differential_gradient}, yields the \emph{power balance equation}
\begin{equation}
    \frac{\d}{\d t}H(\vec{x}(t)) = \nabla H (\vec{x})\transp \dot{\vec{x}} = \vec{z}\transp \vec{E}(\vec{x}) \dot{\vec{x}} = \vec{z}\transp (\vec{J}(\vec{x}) \vec{z} + \vec{B}(\vec{x})\,  \vec{u}) = \vec{y}\transp \vec{u} ,
    \label{eq:power_balance}
\end{equation}
showing passivity and losslessness of the system along every solution $(\vec{x},\vec{u},\vec{y})$ of \eqref{eq:port_hamiltonian_form} for all $t$.
The total change of energy is given by the power supplied through the inputs and outputs.
In the absence of inputs ($\vec{u}=\vec{0}$), the system becomes energy-conserving, i.e., $\dot{H}=0$.

Since energetic consistency is such an important feature of PH systems, much focus in the literature has been devoted to structure-preserving numerical methods that ensure that a discrete counterpart of \eqref{eq:power_balance} holds true also at the discrete level, see e.g. \cite{KotL19,giesselmann_2024_energyconsistent,kinon2025discrete} and references therein.
Throughout this work, we will bring multibody systems into the form \eqref{eq:port_hamiltonian_form} such that their
inherent property of losslessness and passivity is directly encoded
into their representation as PH-DAE. We will then also employ a structure-preserving time discretization scheme.

\section{Port-Hamiltonian formulation of constrained mechanical systems in the Lagrangian setting}\label{sec_PH_lagrangian}




So far, the port-Hamiltonian representation of finite-dimensional mechanical systems largely focused on a definition on the cotangent bundle. In this case, the dynamics are formulated in terms of generalized coordinates $\vec{q}$ and conjugate momenta $\vec{p}$.
To this end, the Legendre transformation to Hamiltonian mechanics is performed and ports are added subsequently.

However, in various works,
the framework of choice is a Lagrangian formulated in terms of generalized coordinates $\vec{q}$ and velocities $\vec{v}$ \cite{betsch2001,betsch_2005_discrete,betsch_2013_consistent,kinon_2023_discrete,kinon2024}.
Therefore, we seek to explore the port-Hamiltonian structure of the Lagrangian equations in a direct fashion. The subsequent treatment is valid for mechanical systems subject to holonomic constraints.
See Figure \ref{fig_flowchart} for an illustration of the two alternative approaches.

\usetikzlibrary{shapes.geometric, arrows}

\tikzstyle{darkblock} = [rectangle, rounded corners,
minimum width=4cm,
minimum height=1cm,
text centered,
draw=black,
text = white,
fill=kit-royalblue70]

\tikzstyle{blueblock} = [rectangle, rounded corners,
minimum width=4cm,
minimum height=1cm,
text centered,
draw=black,
fill=kit-royalblue30]

\tikzstyle{greenblock} = [rectangle, rounded corners,
minimum width=4cm,
minimum height=1cm,
text centered,
draw=black,
fill=kit-green!30]

\tikzstyle{arrow} = [thick,->,>=stealth]

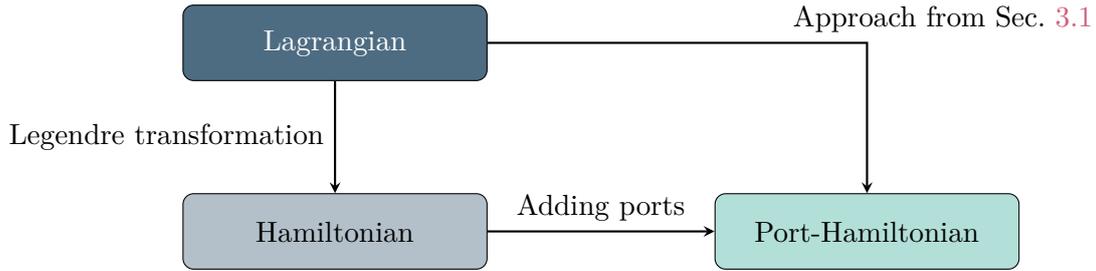
\begin{figure}[t]
    \centering
    \begin{tikzpicture}[node distance=2cm]

        \node (lagrange) [darkblock] {Lagrangian};
        \node (hamilton) [blueblock, below of=lagrange, yshift=-0.5cm] {Hamiltonian};
        \node (phs) [greenblock, right of=hamilton, xshift=5cm] {Port-Hamiltonian};

        \draw [arrow] (lagrange) -- node[anchor=east] {Legendre transformation} (hamilton);
        \draw [arrow] (hamilton) -- node[anchor=south] {Adding ports} (phs);
        \draw [arrow] (lagrange) -| node[right=1cm,anchor=south] {Approach from Sec.~\ref{sec_lagrange}} (phs);

    \end{tikzpicture}
    \caption{Procedure of the present approach. Instead of taking a Hamiltonian approach, we retain the Lagrangian description to directly set up the corresponding port-Hamiltonian formulation.}
    \label{fig_flowchart}
\end{figure}

In Section \ref{sec_lagrange} we introduce a port-Hamiltonian formulation for discrete (possibly constrained) mechanical systems that is equivalent to the Lagrangian equations of the first kind. In Section \ref{sec_robotics} the obtained port-Hamiltonian (or maybe better port-Lagrangian) formulation is specified for a well-known class of multibody systems. The important case of constant mass matrices is highlighted in Section~\ref{sec_const_M}, which will be the basis for the subsequent developments.

\subsection{Constrained mechanical systems in the Lagrangian setting}\label{sec_lagrange}


The motion of conservative mechanical systems with (possibly redundant) position coordinates $\vec{q} \in \mathbb{R}^n$ and velocities $\dot{\vec{q}}=\vec{v} \in \mathbb{R}^n$ is governed by the Lagrange equations of the first kind \cite{simeon_2013_computational,murray_1994_mathematical}.
The stationary conditions are given by the
Euler-Lagrange equations \cite{simeon_2013_computational,murray_1994_mathematical}.
Specifically, taking into account non-conservative forces acting on the system, the equations of motion can be written as
\begin{subequations}
    \label{eq:EL}
    \begin{align}
        \dot{\vec{q}}                                                                  & = \vec{v}     \label{eq:kinematics}                                                                                                                                              \\
        \frac{\d}{\d t}\left(\partial_{\vec{v}} L\left(\vec{q}, \vec{v}\right) \right) & = \partial_{\vec{q}} L \left(\vec{q}, \vec{v}\right) - \nabla \vec{g}(\vec{q})\transp \vec{\lambda} + \tilde{\vec{B}}\left(\vec{q} , \vec{v}\right)\vec{u}  , \label{eq:first_L} \\
        \vec{0}                                                                        & = \nabla \vec{g}(\vec{q}) \vec{v} \label{eq_constraint_velo}
    \end{align}
\end{subequations}
where $L \in \mathcal{C}^1 (\mathbb{R}^n \times  \mathbb{R}^n, \mathbb{R})$ is the Lagrangian function and
$\vec{u} \in \mathbb{R}^p$ are external inputs that are projected using the matrix $\tilde{\vec{B}} \in \mathbb{R}^{n \times p}$ to act as generalized forcing terms (regarded as the known input of the system).
The momentum balance equation contains constraint forces whose magnitude is governed by the Lagrange multipliers $\vec{\lambda} \in \mathbb{R}^m$ which enforce the holonomic and scleronomic constraints
\begin{align} \label{hol_constraint}
    \vec{g}(\vec{q}) = \vec{0} .
\end{align}
Note that the constraints \eqref{hol_constraint} are contained in \eqref{eq_constraint_velo}
in time-differentiated form.
The constraints are assumed to be independent so that the Jacobian $\nabla \vec{g}(\vec{q}) \in \mathbb{R}^{m \times n}$ is of full rank $m$. The constraints
\eqref{eq_constraint_velo} are sometimes called (hidden) velocity-level constraints, as they restrict admissible velocities $\vec{v}$ to the tangent space of the constraint manifold. Equations \eqref{eq:EL} represent differential-algebraic equations with differential index two. For $m=0$, the unconstrained case governed by pure ordinary differential equations is recovered.
Let us now show that the system \eqref{eq:EL} can be represented as a port-Hamiltonian descriptor system in the sense of \eqref{eq:port_hamiltonian_form}.

Applying the chain rule on the left-hand side of \eqref{eq:first_L},
making use of \eqref{eq:kinematics} and adding zero terms $+(\partial^2_{\vec{q}\vec{v}}L) \vec{v}- (\partial^2_{\vec{q}\vec{v}}L) \vec{v}$ on the right-hand side of \eqref{eq:first_L},
\eqref{eq:EL} can be rewritten in the form
\begin{align} \label{eq:new}
    \begin{bmatrix}
        \vec{I} & \vec{0}                       & \vec{0} \\
        \vec{0} & \partial_{\vec{v}\vec{v}}^2 L & \vec{0} \\
        \vec{0} & \vec{0}                       & \vec{0}
    \end{bmatrix} \begin{bmatrix}
                      \dot{\vec{q}} \\ \dot{\vec{v}} \\ \dot{\vec{\lambda}}
                  \end{bmatrix} = \begin{bmatrix}
                                      \vec{0}  & \vec{I}                                                              & \vec{0}                \\
                                      -\vec{I} & \partial^2_{\vec{q}\vec{v}}L - (\partial^2_{\vec{q}\vec{v}}L)\transp & -\nabla \vec{g}\transp \\
                                      \vec{0}  & \nabla \vec{g}                                                       & \vec{0}
                                  \end{bmatrix} \begin{bmatrix}
                                                    -  \partial_{\vec{q}} L + (\partial^2_{\vec{q}\vec{v}}L) \vec{v} \\  \vec{v} \\ \vec{\lambda}
                                                \end{bmatrix}
    + \begin{bmatrix}
          \vec{0} \\ \tilde{\vec{B}}\left(\vec{q} , \vec{v}\right) \\ \vec{0}
      \end{bmatrix} \vec{u} ,
\end{align}
see notation details in Section~\ref{annex_notation}.
Introducing the state vector containing positions, velocities and Lagrange multipliers $\vec{x}=(\vec{q},\vec{v},\vec{\lambda})$, \eqref{eq:new} can be recast as a PH-DAE in the style of \eqref{eq:port_hamiltonian_form}. To see this, we identify
\begin{equation} \label{eq:matrices_E_z}
    \vec{E}(\vec{x}) =
    \begin{bmatrix}
        \vec{I} & \vec{0}                                        & \vec{0} \\
        \vec{0} & \partial_{\vec{v}\vec{v}}^2 L(\vec{q},\vec{v}) & \vec{0} \\
        \vec{0} & \vec{0}                                        & \vec{0}
    \end{bmatrix}
    \quad\mbox{and}\quad
    \vec{z}(\vec{x}) =
    \begin{bmatrix}
        -  \partial_{\vec{q}} L(\vec{q},\vec{v}) + (\partial^2_{\vec{q}\vec{v}}L (\vec{q},\vec{v})) \vec{v} \\  \vec{v} \\ \vec{\lambda}
    \end{bmatrix} .
\end{equation}
The formal definition of a PH system requires the satisfaction of the constitutive relation \eqref{eq:relate_costate_Hamiltonian_differential_gradient} in order to fulfill a power balance equation. Using \eqref{eq:matrices_E_z}, condition \eqref{eq:relate_costate_Hamiltonian_differential_gradient} becomes
\begin{subequations} \label{gradient_hamiltonian}
    \begin{align}
        -\partial_{\vec{q}}L + (\partial_{\vec{q}\vec{v}}^2L) \vec{v} & = \partial_{\vec{q}}H ,       \\
        (\partial_{\vec{v}\vec{v}}^2L) \vec{v}                        & = \partial_{\vec{v}}H ,       \\
        \vec{0}                                                       & = \partial_{\vec{\lambda}}H .
    \end{align}
\end{subequations}
These conditions can be met by choosing the total energy function $H$ as follows:
\begin{equation}\label{eq:legendre}
    H(\vec{q},\vec{v},\vec{\lambda}) = \tilde{H}(\vec{q}, \vec{v}) = \partial_{\vec{v}}{L\left(\vec{q}, \vec{v}\right)}\transp \vec{v}- L\left(\vec{q}, \vec{v}\right) .
\end{equation}
Note that this definition of the total energy $H$ emanates directly from the Lagrangian $L$ and is closely related to the Legendre transformation. We further note that the collocated output, defined through \eqref{eq:relate_costate_output}, is  given by the velocities $\vec{y}=\tilde{\vec{B}}\transp\vec{v}$. Correspondingly, the power balance equation \eqref{eq:power_balance} is given by
\begin{align} \label{eq:power_balance_specific}
    \dot{H} = \vec{u}\transp\tilde{\vec{B}}\transp \vec{v} ,
\end{align}
where the input power results from the scalar product of the external generalized forces $\tilde{\vec{B}}\vec{u}$ with the generalized velocities $\vec{v}$.

\begin{remark}
    One can link the formulation \eqref{eq:new} in the unconstrained case to an equivalent representation,
    where a skew-symmetric matrix appears on the left-hand side of the equation, such that
    \begin{align} \label{eq:winandy}
        \begin{bmatrix}
            \partial_{\vec{v}\vec{q}}^2 L - (\partial_{\vec{v}\vec{q}}^2 L)\transp & \partial_{\vec{v}\vec{v}}^2 L \\
            -(\partial_{\vec{v}\vec{v}}^2 L) \transp                               & \vec{0}
        \end{bmatrix} \begin{bmatrix}
                          \dot{\vec{q}} \\ \dot{\vec{v}}
                      \end{bmatrix} = \begin{bmatrix}
                                          \partial_{\vec{q}} L - (\partial_{\vec{q}\vec{v}}^2 L) \vec{v} \\   -(\partial_{\vec{v}\vec{v}}^2 L) \transp \vec{v}
                                      \end{bmatrix}
        + \begin{bmatrix}
              \tilde{\vec{B}}\left(\vec{q} , \vec{v}\right)\vec{u} \\ \vec{0}
          \end{bmatrix} ,
    \end{align}
    see e.g. \cite[Eq. 1.6]{winandy_dynamics_2019}.
    By comparing \eqref{eq:winandy} with the first two lines of \eqref{gradient_hamiltonian}, the right-hand side of \eqref{eq:winandy} can be identified as the negative gradient of the Hamiltonian \eqref{eq:legendre}. Notably, this system falls into the problem class governed by
    $\vec{D}(\vec{x}) \dot{\vec{x}} = - \nabla H(\vec{x}) + \vec{f}(\vec{x})$,
    which also accommodates gradient systems \cite{EggerHabrichShashkov2021}.
\end{remark}

\subsection{Port-Hamiltonian representation of multibody systems in the Lagrangian setting}\label{sec_robotics}

In the following, we want to specify the newly obtained Lagrangian-based PH-DAE \eqref{eq:new} for the specific class of multibody systems \cite{murray_1994_mathematical,bauchau_2011_flexible,simeon_2013_computational} with the well-known Lagrange function
\begin{align} \label{eq:lagrange_function}
    L(\vec{q}, \vec{v}) = T(\vec{q}, \vec{v}) - V(\vec{q}) ,
\end{align}
where the kinetic energy is given by
\begin{align} \label{eq:kinetic-energy}
    T(\vec{q}, \vec{v}) = \frac{1}{2} \vec{v}\transp \vec{M}(\vec{q}) \vec{v} .
\end{align}

Here, $\vec{M}(\vec{q}) = M_{ij} \vec{e}_i \otimes \vec{e}_j =\vec{M}(\vec{q})\transp \succ 0$ denotes the symmetric and positive definite mass matrix and $V \in \mathcal{C}^1(\mathbb{R}^n, \mathbb{R})$ the potential energy function.
Inserting the Lagrangian \eqref{eq:lagrange_function} into the equations of motion \eqref{eq:EL}, we obtain
\begin{subequations}
    \label{EL_robotics}
    \begin{align}
        \dot{\vec{q}}                 & = \vec{v}   ,                                                                                                                                                                           \\
        \vec{M}(\vec{q})\dot{\vec{v}} & = - \vec{C}(\vec{q}, \vec{v}) \vec{v} - \nabla V(\vec{q}) - \nabla \vec{g}(\vec{q})\transp \vec{\lambda} +\tilde{\vec{B}}\left(\vec{q} , \vec{v}\right) \vec{u} , \label{EL_robotics_2} \\
        \vec{0}                       & =     \nabla \vec{g}(\vec{q}) \vec{v}  \label{EL_robotics_3},
    \end{align}
\end{subequations}
see e.g. \cite[Sect. 7.3]{spong2020robot}, \cite[Sect. 4.5 and 4.6]{vanderschaft2000L2} for the unconstrained case or \cite[Eq. (6.5)]{murray_1994_mathematical} for the constrained case.
In \eqref{EL_robotics}, the term $\vec{C}(\vec{q}, \vec{v})\vec{v}$ denotes the Coriolis and centrifugal forces,
cf. \cite{spong2020robot}.
The corresponding Coriolis matrix $\vec{C}(\vec{q}, \vec{v})$ with its $(k,j)$-th element is given by
\begin{align} \label{coriolis_matrix}
    C_{kj}(\vec{q}, \vec{v}) & = c_{ijk}(\vec{q})v_i .
\end{align}
In \eqref{coriolis_matrix},
we have introduced the
Christoffel symbols of the first kind
\begin{align}
    c_{ijk}(\vec{q}) := \frac{1}{2} \left(\frac{\partial M_{kj}(\vec{q})}{\partial q_i} + \frac{\partial M_{ki}(\vec{q})}{\partial q_j} -  \frac{\partial M_{ij}(\vec{q})}{\partial q_k} \right) ,
\end{align}
satisfying symmetry in the first index pair, i.e. for a fixed $k$ we have $c_{ijk}=c_{jik}$.

The equations of motion \eqref{EL_robotics} can now be expressed equivalently by means of \eqref{eq:new} as a PH-DAE.
From a direct computation of the first and second derivatives of the Lagrangian \eqref{eq:lagrange_function}, formulation \eqref{eq:new} yields
    {\small
        \begin{align} \label{PH_robotics}
            \begin{bmatrix}
                \vec{I} & \vec{0} & \vec{0} \\ \vec{0} & \vec{M}(\vec{q}) & \vec{0} \\ \vec{0} & \vec{0} & \vec{0}
            \end{bmatrix} \begin{bmatrix}
                              \dot{\vec{q}} \\ \dot{\vec{v}} \\ \dot{\vec{\lambda}}
                          \end{bmatrix} = \begin{bmatrix}
                                              \vec{0}  & \vec{I}                                             & \vec{0}                         \\
                                              -\vec{I} & \dot{\vec{M}}(\vec{q}) - 2\vec{C}(\vec{q}, \vec{v}) & -\nabla \vec{g}(\vec{q})\transp \\
                                              \vec{0}  & \nabla \vec{g}(\vec{q})                             & \vec{0}
                                          \end{bmatrix} \begin{bmatrix}
                                                            \partial_{\vec{q}} T(\vec{q}, \vec{v}) + \nabla V(\vec{q}) \\ \vec{v} \\ \vec{\lambda}
                                                        \end{bmatrix} + \begin{bmatrix}
                                                                            \vec{0}                                       \\
                                                                            \tilde{\vec{B}}\left(\vec{q} , \vec{v}\right) \\ \vec{0}
                                                                        \end{bmatrix} \vec{u}.
        \end{align}}
Based on \eqref{eq:lagrange_function}, relation \eqref{eq:legendre} yields the total energy function
\begin{equation}\label{lagrange_storage}
    H\left(\vec{q}, \vec{v}\right) = \frac{1}{2}\vec{v}\transp \vec{M}\left(\vec{q}\right) \vec{v}+ V\left(\vec{q}\right) = T(\vec{q},\vec{v}) + V(\vec{q})
\end{equation}
as the sum of kinetic and potential energy of the system.
The corresponding power balance can now be obtained using two approaches, starting from \eqref{EL_robotics} or from the equivalent \eqref{PH_robotics}.
In both cases, one obtains
\begin{align} \label{eq_power}
    \frac{\d}{\d t} H
    = \vec{v}\transp \tilde{\vec{B}}\vec{u} + \frac{1}{2} \vec{v} \transp \left( \dot{\vec{M}} - 2\vec{C} \right) \vec{v}
    = \vec{v}\transp \tilde{\vec{B}}\vec{u} ,
\end{align}
which is consistent with \eqref{eq:power_balance} and \eqref{eq:power_balance_specific}. The last equality follows from
the skew-symmetry of the term $\dot{\vec{M}}-2\vec{C}$ (see the following Remark or \cite[Proposition 7.1]{spong2020robot}).

\begin{remark}
    The skew-symmetry of the matrix $\dot{\vec{M}}-2\vec{C}$ either follows from the skew-symmetry of the equivalent expression in \eqref{eq:new} or can be verified directly using index notation, such that
    \begin{align}
        \dot{M}_{kj} -2C_{kj} & = \left(\frac{\partial M_{kj}}{\partial q_i} -  \frac{\partial M_{kj}}{\partial q_i} - \frac{\partial M_{ki}}{\partial q_j} + \frac{\partial M_{ij}}{\partial q_k} \right) v_i = \left( \frac{\partial M_{ij}}{\partial q_k} - \frac{\partial M_{ki}}{\partial q_j} \right) v_i \label{skew_symmetry_property} ,
    \end{align}
    where the change of indices $k$ and $j$ flips the sign.
\end{remark}

\subsection{Special case of constant mass matrices}\label{sec_const_M}

In the context of this work it will be sufficient to assume
a constant mass matrix $\vec{M} = \mathrm{const.}$ appearing in the Lagrangian \eqref{eq:lagrange_function},
such that
$ \dot{\vec{M}} = \vec{0}$, $\partial_{\vec{q}}T = \vec{0}$, and $\vec{C}=\vec{0}$.
Correspondingly,  the total energy function \eqref{lagrange_storage}
simplifies to
\begin{align}
    H = \frac{1}{2} \vec{v}\transp \,\vec{M} \,  \vec{v}+ V(\vec{q}) .
    \label{eq:hamiltonian_constrained_mechanical}
\end{align}
and the PH-DAE \eqref{PH_robotics}
takes the simpler form
\begin{align}
    \begin{bmatrix}
        \vec{I} & \vec{0} & \vec{0} \\ \vec{0} & \vec{M} & \vec{0} \\ \vec{0} & \vec{0} & \vec{0}
    \end{bmatrix} \begin{bmatrix}
                      \dot{\vec{q}} \\ \dot{\vec{v}} \\ \dot{\vec{\lambda}}
                  \end{bmatrix} = \begin{bmatrix}
                                      \vec{0}  & \vec{I}                 & \vec{0}                         \\
                                      -\vec{I} & \vec{0}                 & -\nabla \vec{g}(\vec{q})\transp \\
                                      \vec{0}  & \nabla \vec{g}(\vec{q}) & \vec{0}
                                  \end{bmatrix} \begin{bmatrix}
                                                    \nabla V(\vec{q}) \\ \vec{v} \\ \vec{\lambda}
                                                \end{bmatrix} + \begin{bmatrix}
                                                                    \vec{0}                                       \\
                                                                    \tilde{\vec{B}}\left(\vec{q} , \vec{v}\right) \\ \vec{0}
                                                                \end{bmatrix} \vec{u} .
    \label{eq:ph_constraint_compact_simple}
\end{align}
The state-dependency of the structure matrix $\vec{J}(\vec{x})$ reflects the geometric nonlinearity of the PH-DAE.
One can observe that for the unconstrained case, the canonical constant structure matrix is recovered. By comparing this with the general formulation \eqref{eq:port_hamiltonian_form}, one identifies the constant descriptor matrix $\vec{E} = \mathrm{diag}(\vec{I},\vec{M},\vec{0})$ and the
co-state function $\vec{z}=(\nabla V(\vec{q}), \vec{v}, \vec{\lambda})$.
Eventually,
the partial derivatives of the Hamiltonian \eqref{eq:hamiltonian_constrained_mechanical} are given by
\begin{align}
    \frac{\partial H}{\partial \vec{q}} = \nabla V(\vec{q}), \qquad  \qquad
    \frac{\partial H}{\partial \vec{v}} = \vec{M}\vec{v} , \qquad \qquad
    \frac{\partial H}{\partial \vec{\lambda}} =\vec{0},
    \label{eq:co_state_variables}
\end{align}
so that the constitutive relation \eqref{eq:relate_costate_Hamiltonian_differential_gradient} is satisfied.

\section{Director formulation of multibody dynamics in the PH framework}\label{sec:modeling}



In this Section, we introduce a PH-DAE description of  multibody systems, starting from the director formulation of rigid bodies and extending it to  kinematic pairs from which arbitrarily complex multibody systems can be built.

\subsection{Rigid body dynamics based on the director formulation}\label{sec:rigid_body_dynamics}
The rotational dynamics of rigid bodies is a fundamental aspect of multibody systems, and the so-called director formulation \cite{betsch2001,betsch2006} provides a versatile framework to describe the rotational motion of rigid bodies, while
allowing for constant mass matrices as assumed in Section \ref{sec_const_M}.
Our goal is to represent rigid body dynamics within the port-Hamiltonian framework introduced in Section \ref{sec_const_M}. To this end, we recall the essential elements of the director formulation for rigid body dynamics always specifying the related terms in \eqref{eq:ph_constraint_compact_simple}.

\begin{figure}
    \centering
    \def\svgwidth{0.4\textwidth}%
    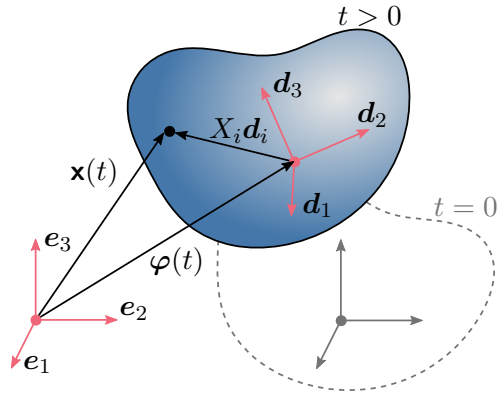%
    \vspace{0.2cm}
    \caption{\small Director formulation of the rigid body \cite{betsch2006}. }
    \label{fig:rigid_body}
\end{figure}

Using the director formulation, the rigid body can be described as a constrained mechanical system, whose equations of motion can be brought into the form of the PH-DAE \eqref{eq:ph_constraint_compact_simple}. The orientation of the rigid body is parameterized by a set of orthonormal director vectors
$\left\{\vec{d}_i\right\}_{i=1}^3$, which are fixed at the center of mass of the rigid body, see Figure \ref{fig:rigid_body}.

For simplicity, the principal axes of inertia are assumed to align with the directors. The directors are subject to
$m=6$ holonomic constraints \eqref{hol_constraint} to ensure the orthogonality and unit length of the directors, i.e.
\begin{align}
    \vec{g}(\vec{q}) = \vec{g}_d(\vec{d}_1, \vec{d}_2, \vec{d}_3) =
    \begin{bmatrix}
        \frac{1}{2} (\vec{d}_1\transp \vec{d}_1 - 1) \\
        \frac{1}{2} (\vec{d}_2\transp \vec{d}_2 - 1) \\
        \frac{1}{2} (\vec{d}_3\transp \vec{d}_3 - 1) \\
        \vec{d}_1\transp \vec{d}_2                   \\
        \vec{d}_1\transp \vec{d}_3                   \\
        \vec{d}_2\transp \vec{d}_3
    \end{bmatrix} = \vec{0} .
    \label{eq:orthonormality_constraints}
\end{align}



With respect to \eqref{eq:ph_constraint_compact_simple}
we are in need of the
gradient of \eqref{eq:orthonormality_constraints}, which is given by
\begin{align}
    \nabla \vec g(\vec q)
    =
    \frac{1}{2}
    \begin{bmatrix}
        \vec 0 & 2 \vec d_1\transp & \vec 0            & \vec 0            \\[0.3em]
        \vec 0 & \vec 0            & 2 \vec d_2\transp & \vec 0            \\[0.3em]
        \vec 0 & \vec 0            & \vec 0            & 2 \vec d_3\transp \\[0.3em]
        \vec 0 & \vec d_2\transp   & \vec d_1\transp   & \vec 0            \\[0.3em]
        \vec 0 & \vec d_3\transp   & \vec 0            & \vec d_1\transp   \\[0.3em]
        \vec 0 & \vec 0            & \vec d_3\transp   & \vec d_2\transp
    \end{bmatrix}.
    \label{eq:orthonormality_constraint_gradient}
\end{align}

The configuration of each body is characterized by the redundant coordinate vector
\begin{align}
    \vec{q}(t) = \begin{bmatrix} \boldsymbol{\varphi}(t) \\ \vec{d}_1(t) \\ \vec{d}_2(t) \\ \vec{d}_3(t) \end{bmatrix} \: \in \: \mathbb{R}^{12},
    \label{eq:director_coordinate_vector}
\end{align}
where $\boldsymbol{\varphi}(t) \in \mathbb{R}^3$ denotes the Cartesian coordinates of the center of mass, formulated with respect to the inertial frame $\{\vec{e}_i\}_{i=1}^3$.
Note that the spatial position of a material point addressed via material coordinates $\vec{X}=X_i \vec{e}_i$ (see Figure \ref{fig:rigid_body}) is given by
\begin{align}
    \vec{\vec{\mathsf{x}}}(\vec{X}, t) = \boldsymbol{\varphi}(t) + \vec{R}(t) \vec{X} \quad \text{with} \quad \vec{R}(t) = \left[ \vec{d}_1(t) \; \vec{d}_2(t) \; \vec{d}_3(t) \right],
    \label{eq:spatial_mapping}
\end{align}


where $\vec{R}(t) \in SO(3)$ is a time-dependent rotation matrix describing the orientation of the director frame, such that
$\vec{d}_i(t) = \vec{R}(t) \, \vec{e}_i$.

We now want to specify the equations of motion for the rigid body formulation. To this end, we can identify the classical relative rotational kinetic energy of a rigid body formulated in convected angular velocities $\vec{\Omega}(t)$ with the kinetic energy expressed in terms of the time derivatives of directors $\dot{\vec{d}}_i(t)$ as
\begin{align} \label{eq:relative-kin-en}
    T_{\mathrm{rot}} = \frac{1}{2} \vec{\Omega}\transp \vec{J} \vec{\Omega} = \frac{1}{2} \sum_{i=1}^3 E_i \, \dot{\vec{d}}_i\transp \dot{\vec{d}}_i ,
\end{align}
where $\vec{J} \in \mathbb{R}^{3 \times 3}$ is the convected inertia tensor of the rigid body.
Additionally, $E_i$ for $i=1,2,3$
are the principal values of the convected Euler tensor \cite{betsch2001} that are related to the principal values of $\vec{J}$ by $E_i = \frac{1}{2}(J_j + J_k - J_i)$ for even permutations of the indices $(i,j,k)$.
Using \eqref{eq:relative-kin-en}, the mass matrix appearing in expression \eqref{eq:hamiltonian_constrained_mechanical} for the total energy can be written as,
\begin{align}
    \vec{M} = \mathrm{diag}(m \vec{I}, E_1 \vec{I}, E_2 \vec{I}, E_3 \vec{I}) .
    \label{eq:mass_matrix}
\end{align}
Notably, this matrix is constant and diagonal, which is a key advantage of the director formulation and can be used for computational efficiency.


External forces and torques acting on the rigid body can be taken into account via the input term $\tilde{\vec{B}}\vec{u}(t)$ in \eqref{eq:ph_constraint_compact_simple}, which contains generalized forces $\vec{f}_{\boldsymbol{\varphi}}(t)$ and $\vec{f}_i(t)$
conjugate to $\boldsymbol{\varphi}$ and $\vec{d}_i$, respectively.

Let $\vec{F}\in\mathbb{R}^3$ be an eccentric force that acts at a distance $\vec{r}=X_i\vec{d}_i$ from the origin of the director frame. In addition to that, let $\vec{\tau}\in\mathbb{R}^3$ be an external torque acting on the rigid body. Then, the total torque relative to the origin of the director frame is given by $\bar{\vec{m}}_\mathrm{ext} = \vec{r}\times \vec{F} + \vec{\tau}$. Now, the generalized forces can be calculated as $\vec{f}_{\boldsymbol{\varphi}} = \vec{F}$ and $\vec{f}_i = -\frac{1}{2} \vec{d}_i\times \bar{\vec{m}}_\mathrm{ext}$, see \cite{betsch2013,betsch2016}. Consequently, arranging $\vec{F}$ and $\vec{\tau}$ in the input vector $\vec{u}=(\vec{F},\vec{\tau}) \in \mathbb{R}^6$, one obtains
\begin{align}
    \begin{bmatrix}
        \vec{f}_{\boldsymbol{\varphi}} \\
        \vec{f}_1                      \\
        \vec{f}_2                      \\
        \vec{f}_3
    \end{bmatrix} = \begin{bmatrix}
                        \vec{F}                                                     \\
                        -\frac{1}{2} \widehat{\vec{d}}_1 \bar{\vec{m}}_\mathrm{ext} \\
                        -\frac{1}{2} \widehat{\vec{d}}_2 \bar{\vec{m}}_\mathrm{ext} \\
                        -\frac{1}{2} \widehat{\vec{d}}_3 \bar{\vec{m}}_\mathrm{ext}
                    \end{bmatrix}
    =
    \begin{bmatrix}
        \vec{I}                                          & \vec{0}                         \\
        -\frac{1}{2}\widehat{\vec{d}}_1\widehat{\vec{r}} & -\frac{1}{2}\widehat{\vec{d}}_1 \\
        -\frac{1}{2}\widehat{\vec{d}}_2\widehat{\vec{r}} & -\frac{1}{2}\widehat{\vec{d}}_2 \\
        -\frac{1}{2}\widehat{\vec{d}}_3\widehat{\vec{r}} & -\frac{1}{2}\widehat{\vec{d}}_3
    \end{bmatrix}
    \begin{bmatrix}
        \vec{F}    \\
        \vec{\tau} \\
    \end{bmatrix} = \tilde{\vec{B}}(\vec{q})\, \vec{u}(t)
    .
    \label{eq:external_force_input_generalized_force}
\end{align}
Therein, $\hat{\square}$ is the mapping \eqref{eq_skew_map} between $\mathbb{R}^3$ and the space of skew-symmetric matrices.
The equations of motion of the free rigid body can now be written in the form of a PH-DAE \eqref{eq:ph_constraint_compact_simple}.

\begin{remark}
    In addition to the power balance equation \eqref{eq:power_balance_specific}, the balance law for the total angular momentum is of crucial importance for rigid body dynamics, see Appendix~\ref{app_ang_mom}. In this connection, the total angular momentum of the rigid body with respect to the origin of the inertial frame is given by
    \begin{align}
        \vec{L} = \boldsymbol{\varphi} \times m \,\dot{\vec{\varphi}} + \sum_{i=1}^{3} \vec{d}_i \times E_i \,\dot{\vec{d}}_i .
        \label{eq:angular_momentum_rigid_body}
    \end{align}
    The related balance equation will be consistently reproduced in discrete time by the structure-preserving integration method presented in \Cref{sec:numerical_integration}.
\end{remark}

\subsection{Multibody systems} \label{sec:MBS}

We extend the rigid body dynamics framework described in the last section to more complex multibody systems comprised of bodies $A, B, ... , N$. The overall system then features
concatenated positions $\vec{q}=(\vec{q}^A, ... , \vec{q}^N)$,
velocities $\vec{v}=(\vec{v}^A, ... , \vec{v}^N)$, and inputs $\vec{u}=(\vec{u}^A, ... , \vec{u}^N)$.
The $N$ rigid bodies in a multibody system are typically connected to each other by joints. In what follows, we focus on kinematic pairs which confine the relative motion of two representative bodies of a multibody system. In the director formulation, the presence of joints can be accounted for by additional constraints.
The constraint vector $\vec{g}$ and the corresponding vector of Lagrange multipliers $\vec{\lambda}$ include contributions from internal director-related constraints, as well as external joint-related constraints, leading to
\begin{align}
    \vec{g} = \begin{bmatrix} \vec{g}_{\text{d}}^{A} \\
                  \vdots                 \\
                  \vec{g}_{\text{d}}^N   \\
                  \vec{g}_{\text{J}}\end{bmatrix}, \qquad
    \nabla \vec{g} = \begin{bmatrix} \nabla \vec{g}_{\text{d}}^{A} \\
                         \vdots                        \\
                         \nabla \vec{g}_{\text{d}}^N   \\
                         \nabla \vec{g}_{\text{J}}\end{bmatrix}, \qquad
    \vec{\lambda} = \begin{bmatrix} \vec{\lambda}_{\text{d}}^{A} \\
                        \vdots                       \\
                        \vec{\lambda}_{\text{d}}^N   \\
                        \vec{\lambda}_{\text{J}}\end{bmatrix} .
    \label{eq:multibody_system_external_constraints_vectors_g_lambda}
\end{align}
In particular, $\vec{g}_{\text{d}}^{k}$ ($k=A,B,\ldots,N$) refers to the internal constraints related to the rigidity of the body $A$, while $\vec{g}_{\text{J}}$ refer to external constraints that arise from the interconnection of bodies in the multibody system through joints.

\begin{figure}
    \centering
    \def\svgwidth{0.7\textwidth}%
    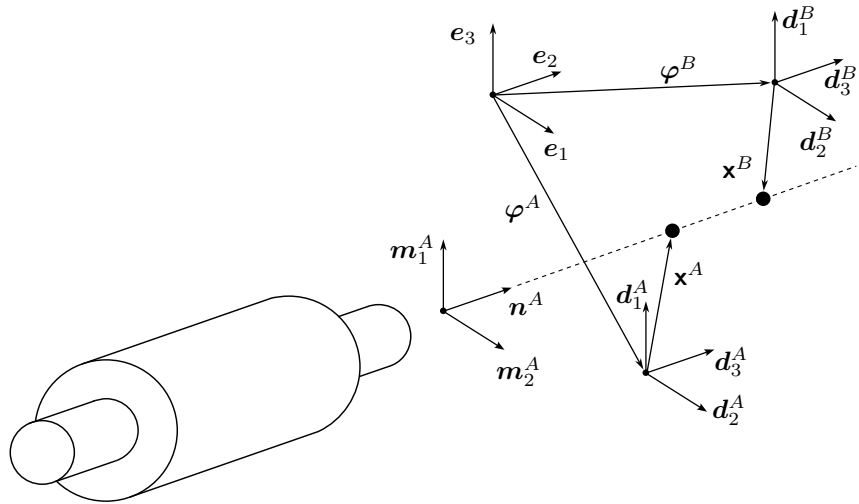%
    \vspace{0.5cm}
    \caption{\small Cylindrical pair in director formulation \cite{betsch2006}. }
    \label{fig:cylindrical_pair}
\end{figure}

As an example, we consider a cylindrical pair connecting two rigid bodies $A$ and $B$ as depicted in Figure \ref{fig:cylindrical_pair}, allowing relative translation along and rotation about a common joint
axis $\vec{n}^A$, which is expressed (without loss of generality) with respect to the coordinates of body $A$. Additional kinematic pairs are summarized in Table~\ref{tab:lower_kinematic_pairs}. For a more in-depth discussion, see \cite{betsch2006,betsch2016}.

\begin{table}[ht]
    \centering
    \caption{Kinematic pairs considered in this work, based on~\cite{betsch2006}.}
    \label{tab:lower_kinematic_pairs}
    \begin{tabular}{lccc}
        \toprule
        Kinematic Pair J & Joint Constraints $m_{\text{J}}$ & DOF of Relative Motion \\
        \midrule
        Spherical (S)    & 3                                & 3                      \\
        Cylindrical (C)  & 4                                & 2                      \\
        Universal (U)    & 4                                & 2                      \\
        Revolute (R)     & 5                                & 1                      \\
        Prismatic (P)    & 5                                & 1                      \\
        \bottomrule
    \end{tabular}
\end{table}

The joint location on a body $\alpha$ is defined by
\begin{align}
    \vec{\mathsf{x}}^\alpha = X_i^\alpha \vec{d}_i^\alpha \qquad \text{for}\qquad \alpha = A, B,
    \label{eq:joint_location}
\end{align}
(see Figure \ref{fig:cylindrical_pair}).
Note that here and in the following, summation over repeated indices is implied unless stated otherwise.
The joint axis vector $\vec{n}^A$ can be regarded as a unit vector with prescribed components $n_i^A$ in the body-fixed frame $\{\vec{d}_i^A\}_{i=1}^3$ such that
\begin{align}
    \vec{n}^A = n_i^A \vec{d}_i^A .
    \label{eq:cylindrical-joint-unit-vector}
\end{align}
Two additional vectors $\vec{m}_\alpha^A = (m_\alpha^A)_i \vec{d}_i^A$ for $\alpha = 1,2$ form a right-handed orthonormal frame given by $\{\vec{n}^A,\vec{m}_1^A,\vec{m}_2^A\}$, see Figure \ref{fig:cylindrical_pair}.
A cylindrical joint gives rise to four external constraint equations
\begin{align}
    \vec{g}_{\text{J,cyl}} (\vec{q}) =
    \begin{bmatrix}
        \vec{m}_1^{A\transp} \Delta \vec{p}     \\
        \vec{m}_2^{A\transp} \Delta \vec{p}     \\
        \vec{n}^{A\transp} \vec{d}_1^B - \eta_1 \\
        \vec{n}^{A\transp} \vec{d}_2^B - \eta_2
    \end{bmatrix} = \vec{0} \quad \text{with} \quad \Delta \vec{p} = \boldsymbol{\varphi}^B + \vec{\mathsf{x}}^B - \boldsymbol{\varphi}^A - \vec{\mathsf{x}}^A,
    \label{eq:cylindrical-pair-constraints}
\end{align}
assuming that $\vec{d}_1^B$ and $\vec{d}_2^B$ are chosen such that they do not align with $\vec{n}^A$.
In \eqref{eq:cylindrical-pair-constraints}, $\Delta \vec{p}$ denotes the relative displacement vector between the two joint locations on bodies $A$ and $B$. The
constants $\eta_1, \eta_2 \in \mathbb{R}$ are chosen to be consistent with initial conditions. While the first two equations in \eqref{eq:cylindrical-pair-constraints} restrict relative translation perpendicular to the joint axis, the last two equations prevent relative rotation about axes orthogonal to the joint axis.

\subsection{Kinematic pairs described as PH interconnection} \label{sec:interconnection}

Let us now gain a different perspective on the modeling of kinematic pairs by employing the notion of power-preserving interconnections within the PH framework \cite{duindam2009,brugnoli2021,berger2025}.
We consider two rigid bodies $A$ and $B$ modeled from the outset as independent PH systems governed by \eqref{eq:port_hamiltonian_form}. That is,
\begin{subequations} \label{eq:port_hamiltonian_form_body_alpha}
    \begin{align}
        \vec{E}^\alpha \dot{\vec{x}}^\alpha & = \vec{J}^\alpha\, \vec{z}^\alpha + \vec{B}^\alpha\,  \vec{u}^\alpha, \label{eq:state_evolution_body_alhpa} \\
        \vec{y}^\alpha                      & = \vec{B}^{\alpha\transp} \vec{z}^\alpha . \label{eq:relate_costate_output_body_alpha}
    \end{align}
\end{subequations}
For each rigid body $\alpha \in \{A, B \}$, the quantities in \eqref{eq:port_hamiltonian_form_body_alpha} have been defined in previous Sections, e.g., in Section \ref{sec_const_M}.
In PH terminology, we call $(\vec{u}^\alpha, \vec{y}^\alpha)$ the port-variables related to interconnection port \cite{vanderschaft2014}, with the associated power $P(\vec{u}^\alpha, \vec{y}^\alpha) = (\vec{u}^\alpha)\transp \vec{y}^\alpha$. This term can be due to external influences (e.g., caused by control inputs) and appears in the power balance \eqref{eq:power_balance}. We show that the PH power-preserving interconnection of these two subsystems, realized via algebraic constraints on the port-variables, is equivalent to the description of kinematic pairs dealt with in the last section.

\begin{figure}
    \centering
    \def\svgwidth{0.5\textwidth}%
\begingroup%
\makeatletter%
\providecommand\color[2][]{%
  \errmessage{(Inkscape) Color is used for the text in Inkscape, but the package 'color.sty' is not loaded}%
  \renewcommand\color[2][]{}%
}%
\providecommand\transparent[1]{%
  \errmessage{(Inkscape) Transparency is used (non-zero) for the text in Inkscape, but the package 'transparent.sty' is not loaded}%
  \renewcommand\transparent[1]{}%
}%
\providecommand\rotatebox[2]{#2}%
\newcommand*\fsize{\dimexpr\f@size pt\relax}%
\newcommand*\lineheight[1]{\fontsize{\fsize}{#1\fsize}\selectfont}%
\ifx\svgwidth\undefined%
  \setlength{\unitlength}{79.73089744bp}%
  \ifx\svgscale\undefined%
    \relax%
  \else%
    \setlength{\unitlength}{\unitlength * \real{\svgscale}}%
  \fi%
\else%
  \setlength{\unitlength}{\svgwidth}%
\fi%
\global\let\svgwidth\undefined%
\global\let\svgscale\undefined%
\makeatother%
\begin{picture}(1,0.55299285)%
  \lineheight{1}%
  \setlength\tabcolsep{0pt}%
  \put(0,0){\includegraphics[width=\unitlength,page=1]{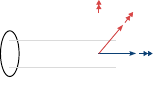}}%
  \put(0.61531494,0.17704134){\color[rgb]{0.04313725,0.25098039,0.45490196}\makebox(0,0)[lt]{\lineheight{1.25}\smash{\begin{tabular}[t]{l} \end{tabular}}}}%
  \put(0.10756845,0.04915219){\color[rgb]{0,0,0}\makebox(0,0)[lt]{\lineheight{1.25}\smash{\begin{tabular}[t]{l}$A$\end{tabular}}}}%
  \put(0.60129563,0.09914959){\color[rgb]{0,0,0}\makebox(0,0)[lt]{\lineheight{1.25}\smash{\begin{tabular}[t]{l}$B$\end{tabular}}}}%
  \put(0,0){\includegraphics[width=\unitlength,page=2]{drawings/kinematic_pair_cylindrical.pdf}}%
  \put(0.10436558,0.12002409){\color[rgb]{0.04705882,0.25098039,0.45490196}\makebox(0,0)[lt]{\lineheight{1.25}\smash{\begin{tabular}[t]{l} \end{tabular}}}}%
  \put(0,0){\includegraphics[width=\unitlength,page=3]{drawings/kinematic_pair_cylindrical.pdf}}%
\end{picture}%
\endgroup%
    \vspace{0.5cm}
    \caption{\small Internal (red) and external (blue) port-variables of a cylindrical pair.}
    \label{fig:cylindrical_pair_ports}
\end{figure}

The interconnection port-variables of the two PH systems at hand are decomposed into \emph{internal} contributions $(\vec{u}_{\text{int}}^\alpha, \vec{y}_{\text{int}}^\alpha)$ and \emph{external} contributions $(\vec{u}_{\text{ext}}^\alpha, \vec{y}_{\text{ext}}^\alpha)$, $\alpha \in \{A, B \}$ such that
\begin{equation} \label{eq:decompose_port}
    \vec{u}^\alpha = \begin{bmatrix}
        \vec{u}_{\text{int}}^\alpha \\
        \vec{u}_{\text{ext}}^{\alpha}
    \end{bmatrix}, \qquad \vec{y}^\alpha = \begin{bmatrix}
        \vec{y}_{\text{int}}^\alpha \\
        \vec{y}_{\text{ext}}^{\alpha}
    \end{bmatrix} .
\end{equation}
The internal contributions are constrained by the interconnection, whereas external contributions remain available for further interaction with the environment or control \cite{berger2025, brugnoli2021}. The distinction becomes more explicit in the example of the cylindrical pair (cf.~Figure \ref{fig:cylindrical_pair_ports}). As described in Section \ref{sec:MBS}, this joint permits body $B$ to translate and rotate along the joint axis only. These admissible motions (marked in blue) correspond to external contributions, since they are not restricted by the interconnection.
Conversely, translations and rotations of body $B$ about the remaining axes are blocked by the joint constraint (marked in red). The associated port-variables therefore constitute the internal contributions.
Body $A$ is not constrained by this interconnection; hence all its associated degrees of freedom appear as external port-variables.

Conforming with \eqref{eq:decompose_port}, we decompose the port matrices in the input-output relations \eqref{eq:port_hamiltonian_form_body_alpha} such that
\begin{align}
    \vec{B}^\alpha \vec{u}^\alpha =
    \begin{bmatrix}
        \vec{B}_{\text{int}}^{\alpha} & \vec{B}_{\text{ext}}^{\alpha}
    \end{bmatrix}
    \begin{bmatrix}
        \vec{u}_{\text{int}}^{\alpha} \\
        \vec{u}_{\text{ext}}^{\alpha}
    \end{bmatrix}, \qquad
    \vec{y}^\alpha =
    \begin{bmatrix}
        \vec{y}_{\text{int}}^{\alpha} \\
        \vec{y}_{\text{ext}}^{\alpha}
    \end{bmatrix} =
    \begin{bmatrix}
        (\vec{B}_{\text{int}}^{\alpha})\transp \\
        (\vec{B}_{\text{ext}}^{\alpha})\transp
    \end{bmatrix}
    \vec{z}^\alpha .
    \label{eq:ph_link_internal_external_ports}
\end{align}
for $\alpha \in \{A, B\}$.
While we can identify $\vec{u}_{\text{int}}^\alpha$ with the forces and torques exerted on body $\alpha$ by the joint, the corresponding $\vec{y}_{\text{int}}^\alpha$ represents the velocity of that body at the joint location.
Now, the interconnection constraints
\begin{subequations} \label{eq:interconnection}
    \begin{align}
        \vec{u}_{\text{int}}^{A} + \vec{u}_{\text{int}}^{B} & = \vec{0}, \label{eq:interconnection_constraint_forces}     \\
        \vec{y}_{\text{int}}^{A} - \vec{y}_{\text{int}}^{B} & = \vec{0}. \label{eq:interconnection_constraint_velocities}
    \end{align}
\end{subequations}
can be applied. Note that \eqref{eq:interconnection_constraint_forces} corresponds to the law of interaction, while \eqref{eq:interconnection_constraint_velocities} confines the relative motion between the two subsystems.
The relations in \eqref{eq:interconnection} correspond to the so-called transformer interconnection \cite{duindam2009,brugnoli2021} between two sets of port-variables given by
\begin{align}
    \begin{bmatrix}
        \vec{y}_\text{int}^{A} \\
        \vec{u}_\text{int}^{A}
    \end{bmatrix}
    = \begin{bmatrix}
          \vec{0}  & \vec{I} \\
          -\vec{I} & \vec{0}
      \end{bmatrix}
    \begin{bmatrix}
        \vec{u}_\text{int}^{B} \\
        \vec{y}_\text{int}^{B}
    \end{bmatrix},
    \label{eq:ch2-transformer-interconnection-id}
\end{align}
showing that the interconnection of the two systems is realized without energy storage or dissipation, so that
\begin{align}
    P_\text{int} = P(\vec{u}_{\text{int}}^{A}, \vec{y}_\text{int}^{A}) + P(\vec{u}_{\text{int}}^{B}, \vec{y}_\text{int}^{B}) = (\vec{u}_\text{int}^{A})\transp \vec{y}_\text{int}^{A} + (\vec{u}_\text{int}^{B})\transp \vec{y}_\text{int}^{B} = 0.
    \label{eq:power-preserving-interconnection}
\end{align}
Making use of the transformer interconnection \eqref{eq:ch2-transformer-interconnection-id} to link the two PH systems of the form \eqref{eq:port_hamiltonian_form_body_alpha}, we
obtain the coupled formulation
\begin{subequations}
    \label{eq:interconnected-system}
    \begin{align}
        \begin{bmatrix}
            \vec{E}^A & \vec{0}   & \vec{0} \\
            \vec{0}   & \vec{E}^B & \vec{0} \\
            \vec{0}   & \vec{0}   & \vec{0}
        \end{bmatrix} \begin{bmatrix}
                          \dot{\vec{x}}^A \\
                          \dot{\vec{x}}^B \\
                          \dot{\vec{\lambda}}_\text{J}
                      \end{bmatrix} & = \begin{bmatrix}
                                            \vec{J}^A                   & \vec{0}                    & \vec{B}_\text{int}^{A}  \\
                                            \vec{0}                     & \vec{J}^B                  & -\vec{B}_\text{int}^{B} \\
                                            -\vec{B}_\text{int}^{A\top} & \vec{B}_\text{int}^{B\top} & \vec{0}
                                        \end{bmatrix}
        \begin{bmatrix}
            \vec{z}^A \\
            \vec{z}^B \\
            \vec{\lambda}_\text{J}
        \end{bmatrix}+ \begin{bmatrix}
                           \vec{B}_\text{ext}^{A} & \vec{0}                \\
                           \vec{0}                & \vec{B}_\text{ext}^{B} \\
                           \vec{0}                & \vec{0}
                       \end{bmatrix} \begin{bmatrix}
                                         \vec{u}_\text{ext}^{A} \\
                                         \vec{u}_\text{ext}^{B}
                                     \end{bmatrix},                 \label{eq:interconnected-system-state-eqs}
        \\
        \begin{bmatrix}
            \vec{y}_\text{ext}^{A} \\
            \vec{y}_\text{ext}^{B}
        \end{bmatrix}              & = \begin{bmatrix}
                                           (\vec{B}_\text{ext}^{A})\transp & \vec{0}                         & \vec{0} \\
                                           \vec{0}                         & (\vec{B}_\text{ext}^{B})\transp & \vec{0}
                                       \end{bmatrix} \begin{bmatrix}
                                                         \vec{z}^A \\
                                                         \vec{z}^B \\
                                                         \vec{\lambda}_\text{J}
                                                     \end{bmatrix}.
    \end{align}
\end{subequations}
With regard to the interconnection constraint \eqref{eq:interconnection_constraint_forces}, we have introduced $\vec{\lambda}_\text{J} := \vec{u}_{\text{int}}^{A} = -\vec{u}_{\text{int}}^{B}$ in the above equation. Furthermore, \eqref{eq:interconnected-system} gives rise to
the total state vector $\vec{x} = (\vec{x}^A, \vec{x}^B, \vec{\lambda}_\text{J})$ and the
Hamiltonian $H(\vec{x}) = H^A(\vec{x}^A) + H^B(\vec{x}^B)$. Note that formulation \eqref{eq:interconnected-system}
fits again into the general framework \eqref{eq:port_hamiltonian_form} and the skew-symmetry of the newly obtained structure matrix encodes \eqref{eq:power-preserving-interconnection}. The remaining power balance still accounts for the external ports only, i.e.
\begin{align}
    \frac{\text{d}}{\text{d}t} H(\vec{x}) = P(\vec{u}_\text{ext}^{A}, \vec{y}_\text{ext}^{A}) + P(\vec{u}_\text{ext}^{B}, \vec{y}_\text{ext}^{B}) =
    \begin{bmatrix}
        \vec{u}_\text{ext}^{A} \\
        \vec{u}_\text{ext}^{B}
    \end{bmatrix}
    \transp
    \begin{bmatrix}
        \vec{y}_\text{ext}^{A} \\
        \vec{y}_\text{ext}^{B}
    \end{bmatrix}
    .
    \label{eq:power-balance-interconnected-system}
\end{align}

Finally, it can be observed that the last row of \eqref{eq:interconnected-system-state-eqs} can be written as
\begin{align}
    \begin{bmatrix}
        -\vec{B}_\text{int}^{A\top} & \vec{B}_\text{int}^{B\top}
    \end{bmatrix}
    \begin{bmatrix}
        \vec{z}^A \\
        \vec{z}^B
    \end{bmatrix} =
    \begin{bmatrix}
        -\tilde{\vec{B}}_\text{int}^{A\transp} & \tilde{\vec{B}}_\text{int}^{B\transp}
    \end{bmatrix} \begin{bmatrix}
                      \vec{v}^A \\
                      \vec{v}^B
                  \end{bmatrix}
    = \vec{0} .
    \label{eq:relation-port-matrix-constraint-gradient}
\end{align}

That is, admissible velocities $\vec{v}=(\vec{v}^A,\vec{v}^B)$ are required to lie in the null space of 
$ \begin{bmatrix} -\tilde{\vec{B}}_\text{int}^{A\transp} & \tilde{\vec{B}}_\text{int}^{B\transp} \end{bmatrix} $. This is in complete analogy to the requirement that the external joint constraints introduced in Section \ref{sec:MBS} be satisfied on the velocity level, which implies 
\begin{equation}
    \frac{\text{d}}{\text{d}t}\vec{g}_{\text{J}} =\nabla \vec{g}_{\text{J}}(\vec{q}) \vec{v} = \vec{0},
\end{equation}
where $\vec{q}=(\vec{q}^A,\vec{q}^B)$. Consequently, 
\begin{equation}
    \text{null} \left(\begin{bmatrix} -\tilde{\vec{B}}_\text{int}^{A\transp} & \tilde{\vec{B}}_\text{int}^{B\transp} \end{bmatrix}\right)
    =
    \text{null} \left(\nabla \vec{g}_{\text{J}}(\vec{q})\right)
\end{equation}


Accordingly, we can conclude that the PH formalism encodes the same information as the conventional formulation of kinematic pairs, but in a more general form that naturally integrates with energy-based interconnection.

\section{Structure-preserving time integration}\label{sec:numerical_integration}
We propose a structure-preserving time-stepping scheme resulting from the discretization of the underlying PH formulation \eqref{eq:port_hamiltonian_form}, e.g. applied to multibody systems \eqref{eq:interconnected-system}. We apply the implicit midpoint rule, 
a second-order accurate integrator. For the present director formulation of multibody systems, this implies the ability of the scheme to preserve consistency with respect to the balance laws for energy and total angular momentum.

Consider the discrete time grid $0 = t^0 < t^1 < \ldots < t^M = t_f$ with $M$ time steps.
Let $t^n$ denote the current time step and define the time step size as $h := t^{n+1} - t^n$.
Further, $\vec{x}^n \approx \vec{x}(t\n)$ and $\vec{x}^{n+1} \approx \vec{x}(t\npe)$ denote the discrete state approximations. The midpoint states are defined by $ \vec{x}^{n+\frac{1}{2}} := \frac{1}{2}(\vec{x}^{n+1} + \vec{x}^n) $.
Applying the midpoint rule to the continuous-time PH equations \eqref{eq:port_hamiltonian_form} yields the time stepping scheme
\begin{subequations} \label{eq:port_hamiltonian_form_discrete}
    \begin{align}
        \vec{E} (\vec{x}^{n+1} - \vec{x}^n) & = h\, \vec{J}(\vec{x}^{n+\frac{1}{2}})\, \vec{z}(\vec{x}^{n+\frac{1}{2}}) + h\, \vec{B}(\vec{x}^{n+\frac{1}{2}})\, \vec{u}^{n+\frac{1}{2}}, \label{eq:time_discrete_phs} \\
        \vec{y}^{n+\frac{1}{2}}             & = \vec{B}(\vec{x}^{n+\frac{1}{2}})\transp \vec{z}(\vec{x}^{n+\frac{1}{2}}), \label{eq:time_discrete_input}
    \end{align}
    such that the discrete-time constitutive relation
    \begin{align}
        \vec{E}\transp \vec{z}(\vec{x}^{n+\frac{1}{2}}) & = \nabla H(\vec{x}^{n+\frac{1}{2}}), \label{eq:time_discrete_costate}
    \end{align}
\end{subequations}
is satisfied by design. Note that we have restricted ourselves here to the case with constant $\vec{E}$. Furthermore, the inputs are evaluated at the midpoint as $\vec{u}^{n+\frac{1}{2}} = \vec{u}(t^{n+\frac{1}{2}})$ but this choice is not mandatory for the structure-preserving properties of the scheme. Correspondingly, the output $\vec{y}^{n+\frac{1}{2}}$ is an approximation of $\vec{y}(t^{n+\frac{1}{2}})$.

In the context of the present PH description of multibody systems, \eqref{eq:time_discrete_phs} represents the discrete counterpart of \eqref{eq:ph_constraint_compact_simple} and
can be given more explicitly as
    {\small
        \begin{align}
            \begin{bmatrix}
                \vec{I} & \vec{0} & \vec{0} \\
                \vec{0} & \vec{M} & \vec{0} \\
                \vec{0} & \vec{0} & \vec{0} \\
            \end{bmatrix}
            \begin{bmatrix}
                \vec{q}\npe - \vec{q}\n \\
                \vec{v}\npe - \vec{v}\n \\
                \vec{\lambda}\npe - \vec{\lambda}\n
            \end{bmatrix}
            = h
            \begin{bmatrix}
                \vec{0}  & \vec{I}                      & \vec{0}                              \\
                -\vec{I} & \vec{0}                      & -\nabla \vec{g}(\vec{q}\npeh)\transp \\
                \vec{0}  & \nabla \vec{g}(\vec{q}\npeh) & \vec{0}
            \end{bmatrix}
            \begin{bmatrix}
                \nabla V(\vec{q}\npeh) \\
                \vec{v}\npeh           \\
                \vec{\lambda}\npeh
            \end{bmatrix}
            + \begin{bmatrix}
                  \vec{0}                       \\
                  \tilde{\vec{B}}(\vec{x}\npeh) \\
                  \vec{0}
              \end{bmatrix} \vec{u}\npeh .
            \label{eq:ph_mbs_compact_discrete}
        \end{align}}

For a quadratic potential function, the scheme \eqref{eq:port_hamiltonian_form_discrete} is a member of the discrete gradient methods proposed in \cite{kinon2025discrete}, yielding exact representations of the power balance equation \eqref{eq:power_balance}.

\subsection{
    Power balance equation in discrete time
}

Note that the structure matrix $\vec{J}(\vec{x})$ in \eqref{eq:port_hamiltonian_form_discrete} is evaluated at the midpoint for second-order convergence, although the following power balance holds independently of the concrete evaluation point.
On the other hand, a midpoint evaluation of the Hamiltonian gradient $\nabla H(\vec{x})$ is crucial for this property, since we assume at most a quadratic Hamiltonian to ensure the directionality property
\begin{align}
    H(\vec{x}^{n+1}) - H(\vec{x}^n) & = \nabla H(\vec{x}^{n+\frac{1}{2}})\transp (\vec{x}^{n+1} - \vec{x}^n) . \label{directionality_property}
\end{align}
Inserting \eqref{eq:port_hamiltonian_form_discrete} results in
\begin{equation}
    \begin{aligned}
        H(\vec{x}^{n+1}) - H(\vec{x}^n) & = h\, (\vec{y}^{n+\frac{1}{2}})\transp \vec{u}^{n+\frac{1}{2}},
    \end{aligned}
    \label{eq:energy_balance_time_discrete}
\end{equation}
where \eqref{eq:time_discrete_costate} and the skew-symmetry of $\vec{J}(\vec{x})$ have been exploited.
This result confirms that the discrete-time formulation preserves the passivity and losslessness of the continuous system. In particular, energy is conserved exactly when $\vec{u} = \vec{0}$.

\subsection{Power-conserving interconnection in discrete time}\label{sec:interconnection_conservation_discrete}

\usetikzlibrary{shapes.geometric, arrows}

\tikzstyle{darkblock} = [rectangle, rounded corners,
minimum width=5cm,
minimum height=1cm,
text centered,
draw=black,
text = white,
fill=kit-royalblue70]

\tikzstyle{blueblock} = [rectangle, rounded corners,
minimum width=5cm,
minimum height=1cm,
text centered,
draw=black,
fill=kit-royalblue30]

\tikzstyle{greenblock} = [rectangle, rounded corners,
minimum width=5cm,
minimum height=1cm,
text centered,
draw=black,
fill=kit-green!30]

\tikzstyle{arrow} = [thick,->,>=stealth]

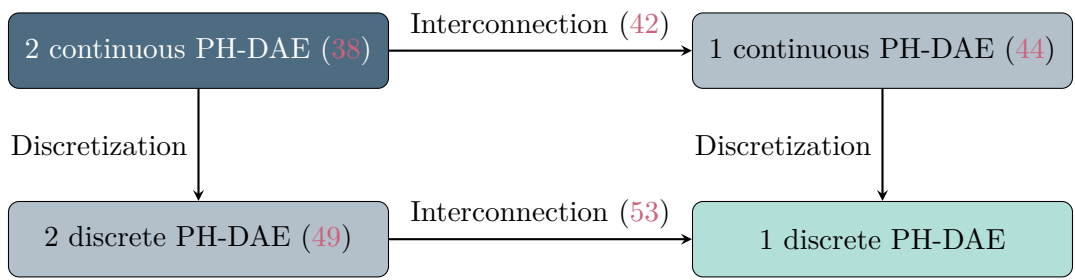
\begin{figure}[b]
    \centering
    \begin{tikzpicture}[node distance=2cm]

        \node (aboveleft) [darkblock] {2 continuous PH-DAE \eqref{eq:port_hamiltonian_form_body_alpha}};
        \node (belowleft) [blueblock, below of=aboveleft, yshift=-0.5cm] {2 discrete PH-DAE \eqref{eq:port_hamiltonian_form_discrete}};
        \node (belowright) [greenblock, right of=belowleft, xshift=7cm] {1 discrete PH-DAE};
        \node (aboveright) [blueblock, right of=aboveleft, xshift=7cm] {1 continuous PH-DAE \eqref{eq:interconnected-system}};

        \draw [arrow] (aboveleft) -- node[anchor=east] {Discretization } (belowleft);
        \draw [arrow] (belowleft) -- node[anchor=south] {Interconnection \eqref{eq:discrete-transformer-interconnection}} (belowright);
        \draw [arrow] (aboveleft) -- node[anchor=south] {Interconnection \eqref{eq:ch2-transformer-interconnection-id}} (aboveright);
        \draw [arrow] (aboveright) -- node[anchor=east]  {Discretization } (belowright);

    \end{tikzpicture}
    \caption{Commuting interconnction and discretization.}
    \label{fig_flowchart_interconnection_discretization}
\end{figure}

The power-conserving nature of the interconnection in continuous time is governed by the transformer relation \eqref{eq:ch2-transformer-interconnection-id}.
It is not obvious how interconnection and time discretization interact with each other.
Therefore, we can show that both operations commute in the sense that interconnecting two subsystems and then applying the time discretization
yields the same result as first discretizing the subsystems and then interconnecting them, see Figure \ref{fig_flowchart_interconnection_discretization}.

This equivalence relies on choosing the same structure-preserving scheme \eqref{eq:port_hamiltonian_form_discrete} for all subsystems.
Consider the discrete formulation \eqref{eq:port_hamiltonian_form_discrete} of two systems with
internal inputs in discrete time $\vec{u}_{\text{int}}^{\alpha,n+\frac{1}{2}}$ and internal outputs $\vec{y}_{\text{int}}^{\alpha,n+\frac{1}{2}}$ for $ \alpha \in \{A, B\} $.
Consider now the interconnection \eqref{eq:ch2-transformer-interconnection-id} in discrete time
\begin{align}
    \begin{bmatrix}
        \vec{y}_\text{int}^{A,n+\frac{1}{2}} \\
        \vec{u}_\text{int}^{A,n+\frac{1}{2}}
    \end{bmatrix}
    =
    \begin{bmatrix}
        \vec{0}  & \vec{I} \\
        -\vec{I} & \vec{0}
    \end{bmatrix}
    \begin{bmatrix}
        \vec{u}_\text{int}^{B,n+\frac{1}{2}} \\
        \vec{y}_\text{int}^{B,n+\frac{1}{2}}
    \end{bmatrix} .
    \label{eq:discrete-transformer-interconnection}
\end{align}
The power exchanged across the interconnection in discrete time is given by
\begin{align}
    (P_\text{int})\npeh = \left(\vec{u}_\text{int}^{A,n+\frac{1}{2}}\right)\transp \vec{y}_\text{int}^{A,n+\frac{1}{2}} + \left(\vec{u}_\text{int}^{B,n+\frac{1}{2}}\right)\transp \vec{y}_\text{int}^{B,n+\frac{1}{2}} =0,
\end{align}
where the discrete interconnection \eqref{eq:discrete-transformer-interconnection} has been used.
This result shows that the discrete-time port-Hamiltonian formulation retains the key power balance property of the continuous system.
Lastly, it can be shown that we obtain the discrete-time system
\begin{align}
    \begin{bmatrix}
        \scriptstyle \vec{E}^A & \scriptstyle \vec{0}   & \scriptstyle \vec{0} \\
        \scriptstyle \vec{0}   & \scriptstyle \vec{E}^B & \scriptstyle \vec{0} \\
        \scriptstyle \vec{0}   & \scriptstyle \vec{0}   & \scriptstyle \vec{0}
    \end{bmatrix} \begin{bmatrix}
                      \scriptstyle \vec{x}^{A,n+1}-\vec{x}^{A,n} \\
                      \scriptstyle \vec{x}^{B,n+1}-\vec{x}^{B,n} \\
                      \scriptstyle \vec{\lambda}_\text{J}^{n+1}-\vec{\lambda}_\text{J}^{n}
                  \end{bmatrix} = \begin{bmatrix}
                                      \scriptstyle \vec{J}^{A,n+\frac{1}{2}}                      & \scriptstyle \vec{0}                                       & \scriptstyle \vec{B}_\text{int}^{A,n+\frac{1}{2}}  \\
                                      \scriptstyle \vec{0}                                        & \scriptstyle \vec{J}^{B,n+\frac{1}{2}}                     & \scriptstyle -\vec{B}_\text{int}^{B,n+\frac{1}{2}} \\
                                      \scriptstyle -(\vec{B}_\text{int}^{A,n+\frac{1}{2}})\transp & \scriptstyle (\vec{B}_\text{int}^{B,n+\frac{1}{2}})\transp & \scriptstyle \vec{0}
                                  \end{bmatrix}
    \begin{bmatrix}
        \scriptstyle h \vec{z}^{A,n+\frac{1}{2}} \\
        \scriptstyle h \vec{z}^{B,n+\frac{1}{2}} \\
        \scriptstyle h \vec{\lambda}_\text{J}^{n+\frac{1}{2}}
    \end{bmatrix}+
    \scriptstyle{h \, \vec{B}_\text{ext}^{n+\frac{1}{2}} \vec{u}_\text{ext}^{n+\frac{1}{2}}} ,
    \label{eq:interconnected-system-discrete}
\end{align}
which is identical to applying the proposed discretization directly to \eqref{eq:interconnected-system-state-eqs}.
The superscript $(\square)^{n+1/2}$ indicates an evaluation of the respective quantity at the midpoint state $\vec{x}\npeh$ and the quantities $\vec{B}_\text{ext}$ and $\vec{u}_\text{ext}$ can be identified by comparison with \eqref{eq:interconnected-system-state-eqs}.

\subsection{Constraint satisfaction in discrete time}\label{sec:constraint_satisfaction_discrete}
We demonstrate that the application of the implicit midpoint rule as outlined above preserves the holonomic constraints exactly in discrete time.
For any at most quadratic constraint function, such as
the director orthonormality constraints \eqref{eq:orthonormality_constraints} or the external joint constraints \eqref{eq:cylindrical-pair-constraints},
it can be shown that in analogy to \eqref{directionality_property} the relation
\begin{align}
    \vec{g}(\vec{q}\npe) - \vec{g}(\vec{q}\n)
    = \nabla \vec{g}(\vec{q}^{n+\frac{1}{2}})(\vec{q}^{n+1}-\vec{q}^n)
\end{align}
holds.
Substituting the kinematic relationship in the first row of \eqref{eq:ph_mbs_compact_discrete} into the above expression yields the expression
\begin{align}
    \vec{g}(\vec{q}\npe) - \vec{g}(\vec{q}\n)
    = h \, \nabla \vec{g}(\vec{q}^{n+\frac{1}{2}}) \vec{v}\npeh = \vec{0} ,
\end{align}
which vanishes due to the third row in \eqref{eq:ph_mbs_compact_discrete}.
Consequently, if the initial values are chosen in a consistent way such that $\vec{g}(\vec{q}^0)=0$,
then the constraints remain satisfied at all discrete times, i.e.
\begin{align}
    \vec{g}(\vec{q}\npe) = \vec{0} , \label{eq_no_drift}
\end{align}
for all $n$.
This includes the orthogonality constraints of the director formulation as well as the joint constraints related to the kinematic pairs contained in Table~\ref{tab:lower_kinematic_pairs}.

\section{Index-reduction by using the GGL principle}\label{sec:GGL_principle}

In view of the exact fulfillment of quadratic holonomic constraints in the sense of \eqref{eq_no_drift}, the index-2 DAE in continuous time \eqref{eq:ph_constraint_compact_simple} become equivalent to a corresponding integration of the related index-3 DAE, see e.g. \cite{kunkel_2006_differentialalgebraic}. It might be beneficial to target a further index reduction also in the PH setting. To this end, we aim for the inclusion of the 
procedure outlined in \cite{kinon2023the-ggl-variational,kinon2023structure-preserving}, where the well-known \emph{Gear-Gupta-Leimkuhler} (GGL) stabilization \cite{gear_1985_automatic}
has been modified variationally.
The equations of motion
resulting from the so-called \emph{GGL variational principle}
represent index-2 DAEs given by
\begin{subequations}
    \begin{align}
        \dot{\boldsymbol{q}}        & = \boldsymbol{v}+\boldsymbol{M}^{-1}\partial_{\boldsymbol{v}}\boldsymbol{g}_v^{\top} \vec{\gamma}     \label{eq_GGL_orig_q} \\
        \vec{M}\dot{\boldsymbol{v}} & = -\nabla V (\vec{q})
        -\nabla \vec{g}(\boldsymbol{q})^{\top}\boldsymbol{\lambda}
        - \partial_{\boldsymbol{q}}\boldsymbol{g}_v^{\top}\vec{\gamma }    \label{eq_GGL_orig_v}                                                                  \\
        \mathbf{0}                  & = \boldsymbol{g}(\boldsymbol{q})  \label{eq_GGL_orig_g}                                                                     \\
        \mathbf{0}                  & = \boldsymbol{g}_v(\boldsymbol{q},\boldsymbol{v}) = \nabla \vec{g}(\boldsymbol{q})\,\boldsymbol{v}. \label{eq_GGL_orig_gv}
    \end{align}
\end{subequations}
These equations extend the well-known index-3 DAEs by the velocity level constraints $\vec{g}_v$, see \eqref{eq_GGL_orig_v}, with additional Lagrange multipliers $\vec{\gamma}$. We omit external generalized forces here for the sake of conciseness.
In this context, the classical GGL--stabilization is modified by one additional term in the balance of momentum to obtain a Hamiltonian and thus symplectic structure \cite{kinon2023structure-preserving}. To fit this framework into the PH structure \eqref{eq:port_hamiltonian_form}, we  differentiate both position and velocity level constraints \eqref{eq_GGL_orig_g} and \eqref{eq_GGL_orig_gv} in time to obtain
\begin{subequations}
    \begin{align}
        \frac{\d}{\d t}\vec{g}(\vec{q})           & = \nabla \vec{g}(\vec{q})\left(  \boldsymbol{v}+\boldsymbol{M}^{-1}\partial_{\boldsymbol{v}}\boldsymbol{g}_v^{\top} \vec{\gamma} \right) = \vec{0}    ,                                                              \\
        \frac{\d}{\d t}\vec{g}_v(\vec{q},\vec{v}) & = \partial_{\vec{q}}\vec{g}_v \left(  \boldsymbol{v}+\boldsymbol{M}^{-1}\partial_{\boldsymbol{v}}\boldsymbol{g}_v^{\top} \vec{\gamma} \right) + \partial_{\vec{v}}\vec{g}_v \vec{M}^{-1} \left(  -\nabla V (\vec{q})
        -\nabla \vec{g}(\boldsymbol{q})^{\top}\boldsymbol{\lambda}
        - \partial_{\boldsymbol{q}}\boldsymbol{g}_v^{\top}\vec{\gamma}  \right) = \vec{0} ,
    \end{align}
\end{subequations}
where we have substituted \eqref{eq_GGL_orig_q} and \eqref{eq_GGL_orig_v}. Eventually, these constraints on velocity and acceleration level are
embedded in the index-1 PH-DAE
\begin{align}
    \begin{bmatrix}
        \vec{I} & \vec{0} & \vec{0} & \vec{0} \\
        \vec{0} & \vec{M} & \vec{0} & \vec{0} \\
        \vec{0} & \vec{0} & \vec{0} & \vec{0} \\
        \vec{0} & \vec{0} & \vec{0} & \vec{0}
    \end{bmatrix}
    \begin{bmatrix}
        \dot{\vec{q}} \\\dot{\vec{v}}\\\dot{\vec{\lambda}}\\\dot{\vec{\gamma}}
    \end{bmatrix}=
    \begin{bmatrix}
        \vec{0}                      & \vec{I}                       & \vec{0}                                               & \vec{M}^{-1} \nabla \vec{g}\transp                \\
        -\vec{I}                     & \vec{0}                       & -\nabla \vec{g}\transp                                & -\partial_{\vec{q}}\,\vec{g}_v\transp             \\
        \vec{0}                      & \nabla \vec{g}                & \vec{0}                                               & \nabla \vec{g} \vec{M}^{-1} \nabla \vec{g}\transp \\
        - \nabla \vec{g}\vec{M}^{-1} & \partial_{\vec{q}}\,\vec{g}_v & -\nabla \vec{g}\vec{M}^{-1}  \: \nabla \vec{g}\transp & \vec{J}_{44}
    \end{bmatrix}
    \begin{bmatrix}
        \nabla V      \\
        \vec{v}       \\
        \vec{\lambda} \\
        \vec{\gamma}
    \end{bmatrix}
    \label{eq:index-1-PH-DAEs-velocity}
\end{align}

where $\vec{J}_{44}(\vec{q},\vec{v}) = \partial_{\vec{q}}\,\vec{g}_v \:\vec{M}^{-\top} \nabla \vec{g}\transp -\nabla \vec{g}\vec{M}^{-1} \:\partial_{\vec{q}}\,\vec{g}_v\transp  = - \vec{J}_{44}\transp$ and $\partial_{\vec{v}} \vec{g}_v = \nabla \vec{g}$. Note that \eqref{eq:index-1-PH-DAEs-velocity} can be viewed as extension of \eqref{eq:ph_constraint_compact_simple}. As these equations again fit into the PH framework \eqref{eq:port_hamiltonian_form}, i.e. $\vec{E}\dot{\vec{x}} = \vec{J}(\vec{x})\vec{z}(\vec{x})$, we can be sure that the power balance equation \eqref{eq:power_balance} remains satisfied.
Therefore, we apply a time integration approach as in \eqref{eq:port_hamiltonian_form_discrete}, using a midpoint evaluation of the structure matrix such that $\vec{E}(\vec{x}\npe - \vec{x}\n) = h \vec{J}(\vec{x}\npeh)\vec{z}(\vec{x}\npeh)$. Analogously to the previous case, we can verify that this formulation is now able to preserve both position and velocity level constraints exactly, i.e. for all $n$
\begin{align} \label{ggl_discrete_satisfaction}
    \vec{g}(\vec{q}\npe) = \vec{0} , \qquad \text{and} \qquad \vec{g}_v(\vec{q}\npe,\vec{v}\npe) = \vec{0} ,
\end{align}
for consistent initial conditions $(\vec{q}^0,\vec{v}^0)$. It is also required that the holonomic constraints are polynomials of order of at most two. The proof follows the lines in \Cref{sec:constraint_satisfaction_discrete}, by including $\vec{g}_v$ in addition to $\vec{g}$.

\section{Numerical results and discussion}\label{sec:numerical_results}

In the following, we apply the proposed PH midpoint integrator \eqref{eq:ph_mbs_compact_discrete}, henceforth labeled as \enquote{PH-MP}, to three representative numerical examples dealing with multibody systems.
Further, we highlight stability benefits that can be gained by extending the integrator to include the GGL-type stabilization introduced in Section \ref{sec:GGL_principle}. We label this integrator as \enquote{PH-MP-GGL}.
Section~\ref{sec:cylindrical_pair} demonstrates the operation of a cylindrical joint as introduced in Section~\ref{sec:interconnection}.
Section~\ref{sec:closed_loop_MBS} discusses a model problem involving a closed-loop system.
Closed-loop systems, in particular, make it challenging to determine minimal coordinates, thereby illustrating the advantages of using redundant coordinates.
Section~\ref{sec:slider_crank_mechanism} investigates the spatial slider-crank mechanism -- a more complex closed-loop system, including gravity and fixed supports.

Newton’s method is used to solve the resulting nonlinear equations in each time step with an absolute tolerance of $\epsilon = 10^{-9}$.
Time-step sizes follow values from the literature unless otherwise stated, in which case they are chosen to ensure convergence.
Notably, for the midpoint integrators at hand, the step size is constrained only by the convergence of Newton’s method at each iteration.

\subsection{Flying cylindrical pair}\label{sec:cylindrical_pair}
\begin{figure}[H]
    \centering
    \begin{subfigure}[t]{0.46\textwidth}
        \centering
        \def\svgwidth{\textwidth}%
        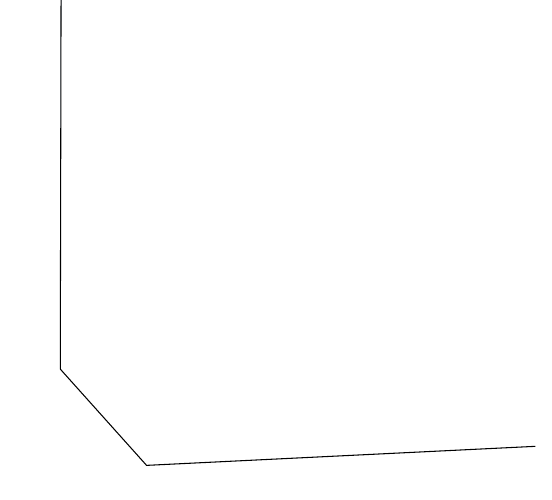%
        \caption{Initial configuration as in \cite{betsch2006}.}
        \label{fig:cylindrical_pair_t0}
    \end{subfigure}
    \hfill
    \begin{subfigure}[t]{0.46\textwidth}
        \centering
        \def\svgwidth{\textwidth}%
        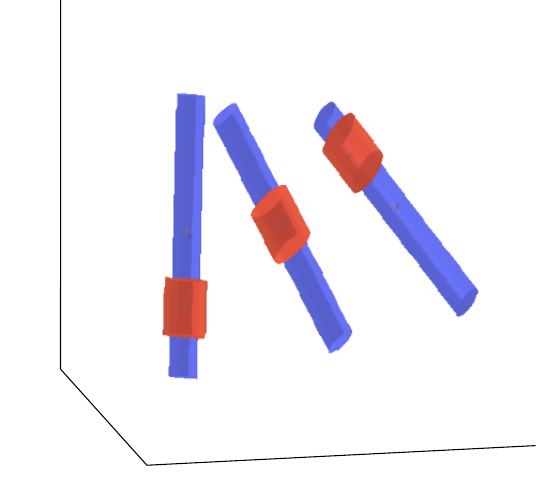%
        \caption{Snapshots of the motion}
        \label{fig:cylindrical_pair_motion}
    \end{subfigure}
    \caption{Flying cylindrical pair.}\label{fig:cylindrical_pair_sketch}
\end{figure}

This example, adapted from \cite{betsch2006}, demonstrates the dynamics of a cylindrical pair in free flight and serves to validate the conservation properties of the proposed framework.
The set-up allows assessing the accuracy of the proposed numerical integration in preserving energy and total angular momentum in the absence of external forces, see \eqref{eq:energy_balance_time_discrete} and \eqref{discrete_angmom_balance}, respectively.
\begin{table}[H]
    \centering
    \caption{Simulation parameters flying cylindrical pair.}
    \label{tab:cylindrical_pair_setup}
    \renewcommand{\arraystretch}{1.3}
    \begin{tabular}{cccccccccc}
        \hline
        $h$   & $t_\text{end}$ & $l^A$ & $r^A$ & $m^A$ & $l^B$ & $r_\text{outer}^B$ & $r_\text{inner}^B$ & $m^B$ \\
        \hline
        0.001 & 0.7            & 30    & 2     & 4     & 6     & 3                  & 2                  & 3     \\
        \hline
    \end{tabular}
\end{table}
The initial configuration of the system is depicted in Figure \ref{fig:cylindrical_pair_t0} and chosen simulation parameters are summarized in Table \ref{tab:cylindrical_pair_setup}.
Body $A$ is a long, slender cylinder with mass $m^A$, length $l^A$ and radius $r^A$, while body $B$ is a short, thick hollow cylinder with mass $m^B$, length $l^B$, inner radius $r_\text{inner}^B$ and outer radius $r_\text{outer}^B$. The cylindrical joint located at the center of mass of body $A$ connects the two bodies, allowing relative rotation about
and translation along the joint axis $\vec{n}^A$ for $\vec{n}^A(t=0)=[
    0 \ 0 \ 1]\transp$. The independent generalized velocities at $t=0$ for the problem based on \cite{betsch2006} are given as
\begin{align*}
    \vec{v}_\text{ind}(t=0) & = \begin{bmatrix}
                                    \left(\vec{v}_{\varphi}^A(t=0)\right)\transp &
                                    \left(\vec{\omega}^A(t=0)\right)\transp      &
                                    \dot{u}^B (t=0)                              &
                                    \dot{\theta}^B(t=0)
                                \end{bmatrix}\transp \\
                            & = \begin{bmatrix}
                                    0    &
                                    50   &
                                    0    &
                                    1    &
                                    1.5  &
                                    0    &
                                    35.5 &
                                    -100
                                \end{bmatrix}\transp.
\end{align*}
Here, $\dot{u}^B$ indicates the translation velocity of body $B$ along $\vec{n}^A$, and $\dot{\theta}^B$ accounts for the angular velocity of body $B$ relative to body $A$. Based on \cite{betsch2006}, we obtain the initial velocity of body $B$ in director formulation using
\begin{align*}
    \vec{v}_{\varphi}^B = \vec{v}_{\varphi}^A+\vec{\omega}^A \times (\vec{\varphi}^B-\vec{\varphi}^A)+\dot{u}^B \vec{n}^A+\dot{\theta}^B \vec{\mathsf{x}}^B\times \vec{n}^A.
\end{align*}
No external forces are acting on the bodies, and no gravitational forces are active.  
Snapshots in Figure~\ref{fig:cylindrical_pair_motion} show the motion of the two bodies at $t \in \left\{0.1,\, 0.4, 0.7\right\}$. With no external forces or gravity acting on the system, both the total energy and
the angular momentum vector $\vec{L}$ remain constant down to machine precision throughout the simulation, see Figures \ref{fig:energy_cylindrical_pair} and \ref{fig:angular_momentum_cylindrical_pair}, respectively.

\begin{figure}[H]
    \centering
    \begin{subfigure}[t]{0.48\textwidth}
        \centering
        \adjustbox{valign=t}{%
            \definecolor{mycolor1}{rgb}{0.90196,0.62353,0.00000}%
\definecolor{mycolor2}{rgb}{0.33725,0.70588,0.91373}%
\definecolor{mycolor3}{rgb}{0.80000,0.47451,0.65490}%

\begin{tikzpicture}
    \begin{axis}[%
            width=0.65\textwidth,
            height=0.1\textheight,
            scale only axis,
            xmin=0,
            xmax=0.7,
            xtick ={0,0.7},
            unbounded coords=discard,
            xlabel={$t^n$},
            ylabel={$H(\vec{x}\n)$},
            ymin = 90000,
            ymax = 120000,
            xmajorgrids,
            ymajorgrids,
            axis background/.style={fill=white},
        ]

        \addplot [color=kit-royalblue100, line width=1.5pt, solid]
        table[col sep=comma, x=time, y=system_total_energy_current_time] {data/e109_cylindrical_pair_Betsch_results_h=0.001_new.csv};

    \end{axis}
\end{tikzpicture}
        }
        \label{fig:total_energy_cylindrical_pair}
    \end{subfigure}
    \hfill
    \begin{subfigure}[t]{0.48\textwidth}
        \centering
        \adjustbox{valign=t}{%
            \definecolor{mycolor1}{rgb}{0.90196,0.62353,0.00000}%
\definecolor{mycolor2}{rgb}{0.33725,0.70588,0.91373}%
\definecolor{mycolor3}{rgb}{0.80000,0.47451,0.65490}%

\begin{tikzpicture}
    \begin{axis}[%
            width=0.65\textwidth,
            height=0.1\textheight,
            scale only axis,
            xmin=0,
            xmax=0.7,
            xtick ={0,0.7},
            xlabel={$t^n$},
            ylabel={$H(\bm{x}^{n+1})-H(\bm{x}^{n})$},
            ymin = -2e-9,
            ymax = 2e-9,
            ytick = {-2e-9, -1e-9, 0, 1e-9, 2e-9},
            xmajorgrids,
            ymajorgrids,
            axis background/.style={fill=white},
        ]

        \addplot [color=kit-royalblue100, line width=1.0pt, solid]
        table[col sep=comma, x=time, y=system_total_energy_interval_increment,restrict y to domain=-1:1] {data/e109_cylindrical_pair_Betsch_results_h=0.001_new.csv};

    \end{axis}
\end{tikzpicture}
        }
        \label{fig:energy_increment_cylindrical_pair}
    \end{subfigure}
    \caption{Hamiltonian for flying cylindrical pair  (using PH-MP).}\label{fig:energy_cylindrical_pair}
\end{figure}
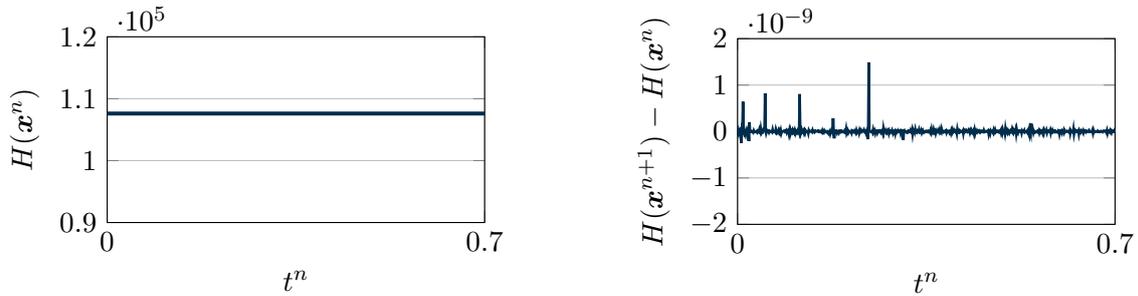
\begin{figure}[H]
    \begin{subfigure}[b]{0.48\textwidth}
        \centering
        \adjustbox{valign=b}{%
            \hspace{-1em}
            \definecolor{mycolor1}{rgb}{0.90196,0.62353,0.00000}%
\definecolor{mycolor2}{rgb}{0.33725,0.70588,0.91373}%
\definecolor{mycolor3}{rgb}{0.80000,0.47451,0.65490}%

\begin{tikzpicture}
    \begin{axis}[%
            width=0.65\textwidth,
            height=0.1\textheight,
            scale only axis,
            ymin=-1200,
            ymax=2200,
            xmin=0,
            xmax=0.7,
            scaled ticks = true,
            xtick ={0,0.7},
            unbounded coords=discard,
            xlabel={$t^n$},
            ylabel={$\vec{L}(\vec{x}\n)$},
            xmajorgrids,
            ymajorgrids,
            axis background/.style={fill=white},
            legend style={at={(0.97,0.03)}, anchor=south east, draw=white!15!black, fill=white, font=\footnotesize},
            legend cell align=left,
            legend entries={$L_x$, $L_y$, $L_z$}
        ]

        \addplot [color=kit-green100, line width=1.5pt, solid, smooth]
        table[col sep=comma, x=time, y=total_angular_momentum_x_current_time] {data/e109_cylindrical_pair_Betsch_results_h=0.001_new.csv};

        \addplot [color=kit-royalblue100, line width=1.5pt, solid, smooth]
        table[col sep=comma, x=time, y=total_angular_momentum_y_current_time] {data/e109_cylindrical_pair_Betsch_results_h=0.001_new.csv};

        \addplot [color=kit-royalblue50, line width=1.5pt, solid, smooth]
        table[col sep=comma, x=time, y=total_angular_momentum_z_current_time] {data/e109_cylindrical_pair_Betsch_results_h=0.001_new.csv};

    \end{axis}
\end{tikzpicture}
        }
        \label{fig:angular_momentum_conservation_cylindrical_pair}
    \end{subfigure}
    \hfill
    \begin{subfigure}[b]{0.48\textwidth}
        \centering
        \adjustbox{valign=b}{%
            \definecolor{mycolor1}{rgb}{0.90196,0.62353,0.00000}%
\definecolor{mycolor2}{rgb}{0.33725,0.70588,0.91373}%
\definecolor{mycolor3}{rgb}{0.80000,0.47451,0.65490}%

\begin{tikzpicture}
    \begin{axis}[%
            width=0.65\textwidth,
            height=0.1\textheight,
            scale only axis,
            xmin=0,
            xmax=0.7,
            ymin = -2e-11,
            ymax = 2e-11,
            ytick ={-2e-11,-1e-11,0,1e-11,2e-11},
            xtick ={0,0.7},
            unbounded coords=discard,
            xlabel={$t^n$},
            ylabel={$\bm{L}(\bm{x}^{n+1})-\bm{L}(\bm{x}^{n})$},
            xmajorgrids,
            ymajorgrids,
            axis background/.style={fill=white},
            legend style={at={(0.97,0.03)}, anchor=south east, draw=white!15!black, fill=white, font=\footnotesize},
            legend cell align=left,
            legend entries={$L_x$, $L_y$, $L_z$}
        ]

        \addplot [color=kit-green100, line width=1.0pt, solid, smooth]
        table[col sep=comma, x=time, y=total_angular_momentum_x_interval_increment] {data/e109_cylindrical_pair_Betsch_results_h=0.001_new.csv};

        \addplot [color=kit-royalblue100, line width=1.0pt, solid, smooth]
        table[col sep=comma, x=time, y=total_angular_momentum_y_interval_increment] {data/e109_cylindrical_pair_Betsch_results_h=0.001_new.csv};

        \addplot [color=kit-royalblue50, line width=1.0pt, solid, smooth]
        table[col sep=comma, x=time, y=total_angular_momentum_z_interval_increment] {data/e109_cylindrical_pair_Betsch_results_h=0.001_new.csv};

    \end{axis}
\end{tikzpicture}
        }
        \label{fig:angular_momentum_increment_cylindrical_pair}
    \end{subfigure}
    \caption{Total angular momentum for cylindrical pair (using PH-MP).}
    \label{fig:angular_momentum_cylindrical_pair}
\end{figure}
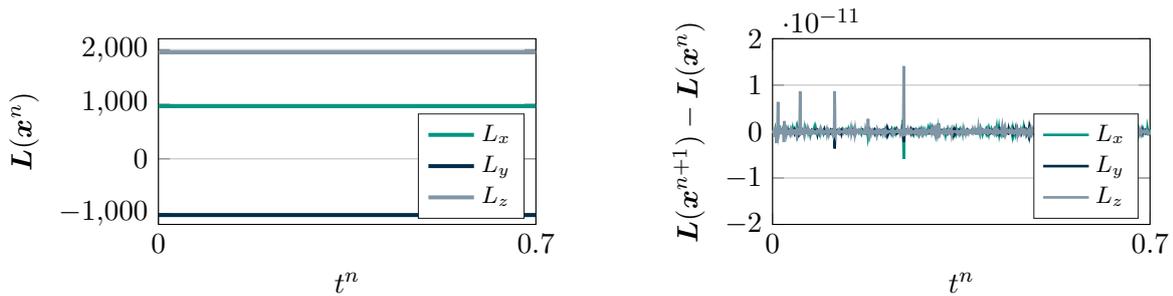
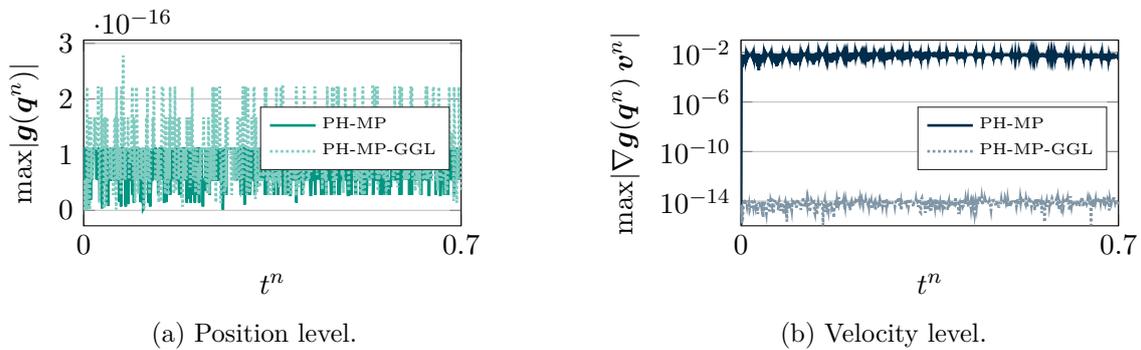
\begin{figure}[H]
    \centering
    \begin{subfigure}[b]{0.48\textwidth}
        \centering
        \adjustbox{valign=b}{%
            \definecolor{mycolor1}{rgb}{0.90196,0.62353,0.00000}%
\definecolor{mycolor2}{rgb}{0.33725,0.70588,0.91373}%
\definecolor{mycolor3}{rgb}{0.80000,0.47451,0.65490}%

\begin{tikzpicture}
    \begin{axis}[%
            width=0.65\textwidth,
            height=0.1\textheight,
            scale only axis,
            xmin=0,
            xmax=0.7,
            xtick ={0,0.7},
            unbounded coords=discard,
            xlabel={$t^n$},
            ylabel={$\text{max} \lvert \bm{g}(\vec{q}\n) \rvert$},
            legend style={at={(0.97,0.3)}, anchor=south east, draw=white!15!black, fill=white, font=\tiny},
            legend cell align=left,
            legend entries={PH-MP, PH-MP-GGL},
            xmajorgrids,
            ymajorgrids,
            axis background/.style={fill=white},
        ]

        \addplot [color=kit-green100, line width=1.0pt, solid]
        table[col sep=comma, x=time, y=max_constraint_error_q] {data/e109_cylindrical_pair_Betsch_results_h=0.001_new.csv};

        \addplot [color=kit-green50, line width=1.0pt, densely dotted]
        table[col sep=semicolon, x=time, y=max_constraint_error_q] {data/e104.2_GGL_cylindrical_pair_Betsch_midpoint_results_t=0.7_h=0.001.csv};

    \end{axis}
\end{tikzpicture}
        }
        \caption{Position level.}
        \label{fig:q_constraint_cylindrical_pair}
    \end{subfigure}
    \hfill
    \begin{subfigure}[b]{0.48\textwidth}
        \centering
        \adjustbox{valign=b}{%
            \definecolor{mycolor1}{rgb}{0.90196,0.62353,0.00000}%
\definecolor{mycolor2}{rgb}{0.33725,0.70588,0.91373}%
\definecolor{mycolor3}{rgb}{0.80000,0.47451,0.65490}%

\begin{tikzpicture}
    \begin{axis}[%
        width=0.65\textwidth,
        height=0.1\textheight,
        scale only axis,
        xmin=0,
        xmax=0.7,
        xtick={0,0.7},
        ymode=log,
        ymode=log,
        ymin=1e-16,
        ymax=1e-1,
        ytickten={-14,-10,-6,-2},
        unbounded coords=discard,
        xlabel={$t^n$},
        ylabel={$\max \lvert \nabla \vec{g}(\vec{q}^n)\,\vec{v}^n\rvert$},
        legend style={at={(0.97,0.3)}, anchor=south east, draw=white!15!black, fill=white, font=\tiny},
        legend cell align=left,
        legend entries={PH-MP, PH-MP-GGL},
        xmajorgrids,
        ymajorgrids,
        axis background/.style={fill=white},
    ]

        \addplot [color=kit-royalblue100, line width=1.0pt, solid]
        table[col sep=comma, x=time, y=max_constraint_error_v] {data/e109_cylindrical_pair_Betsch_results_h=0.001_new.csv};

        \addplot [color=kit-royalblue50, line width=1.0pt, densely dotted]
        table[col sep=semicolon, x index=80, y index=146] {data/e104.2_GGL_cylindrical_pair_Betsch_midpoint_results_t=0.7_h=0.001.csv};

    \end{axis}
\end{tikzpicture}
        }
        \caption{Velocity level.}
        \label{fig:v_constraint_cylindrical_pair}
    \end{subfigure}
    \caption{Constraint satisfaction (maximum absolute error) for cylindrical pair.}
    \label{fig:constraint_satisfaction_cylindrical_pair}
\end{figure}

Analyzing the constraint satisfaction, as discussed in \eqref{eq_no_drift}, Figure \ref{fig:q_constraint_cylindrical_pair} shows that position-level constraints remain satisfied up to numerical precision throughout the simulation,
confirming the effectiveness of the proposed approach in maintaining constraint consistency, see Section \ref{sec:constraint_satisfaction_discrete}. As depicted in Figure \ref{fig:v_constraint_cylindrical_pair},
this is not the case for constraints on velocity level with the PH-MP integrator. In contrast, the PH-MP-GGL integrator exactly enforces constraints also on velocity-level, see \eqref{ggl_discrete_satisfaction}, ensuring that they are satisfied up to numerical precision, cf. Figure \ref{fig:v_constraint_cylindrical_pair}.
In order to verify the accuracy of the numerical scheme, the convergence of the state variables $\vec{q},\,\vec{v},\,\vec{\lambda}$ and the conserved quantities $H,\,\vec{L}$ is analyzed.
Results for $h=\left\{10^{-2}, 10^{-3}, 10^{-4}\right\}$ are compared to a finely discretized solution with $h=10^{-5}$. For each state, we calculate the root-mean-square (RMS) error between the simulation results and
the reference solution at $t=\bar{t}=0.02$ as, e.g.,
\begin{align}
    e_\text{RMS}(\vec{q}) = \sqrt{\frac{1}{n}\sum_{i=1}^{N}\left(\vec{q}_i^{\text{sim}}(\bar{t}) - \vec{q}_i^{\text{ref}}(\bar{t})\right)^2}.
    \label{eq:rms_q}
\end{align}
The calculation of the RMS error for velocity $\vec{v}$, Lagrange multiplier $\vec{\lambda}$ and the conserved quantities $H,\,\vec{L}$ follows the same procedure.

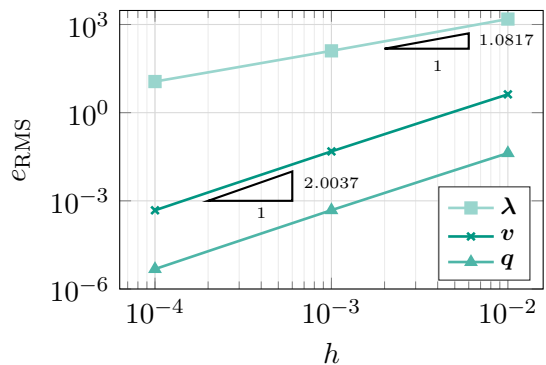
\begin{figure}[H]
    \centering
    \centering
    \adjustbox{valign=b}{%
        \definecolor{mycolor1}{rgb}{0.90196,0.62353,0.00000}%
\definecolor{mycolor2}{rgb}{0.33725,0.70588,0.91373}%
\definecolor{mycolor3}{rgb}{0.80000,0.47451,0.65490}%

\begin{tikzpicture}
    \begin{loglogaxis}[%
            width=0.35\textwidth,
            height=0.15\textheight,
            scale only axis,
            ymin=0.000001,
            ymax=3000,
            unbounded coords=discard,
            xlabel={$h$},
            ylabel={$e_{\text{RMS}}$},
            grid=both,                   
            major grid style={gray!30},  
            minor grid style={gray!15},  
            axis background/.style={fill=white},
            legend style={at={(0.97,0.03)}, anchor=south east, draw=white!15!black, fill=white, font=\footnotesize},
            legend cell align=left,
            legend entries={$\bm{\lambda}$, $\bm{v}$, $\bm{q}$}
        ]

        \addplot [mark=square*,  mark size=2pt, color=kit-green40, line width=1pt, solid, smooth]
        table[col sep=comma,trim cells, x index=0, y index=5] {data/convergence_cylindrical_pair.csv};

        \addplot [mark=x, mark size=2pt,color=kit-green100, line width=1pt, solid, smooth]
        table[col sep=comma,trim cells, x index=0, y index=2] {data/convergence_cylindrical_pair.csv};

        \addplot [mark=triangle*,  mark size=2pt, color=kit-green70, line width=1pt, solid, smooth]
        table[col sep=comma,trim cells, x index=0, y index=1] {data/convergence_cylindrical_pair.csv};

        \coordinate (A) at (axis cs:2e-4,1e-3);
        \coordinate (B) at (axis cs:6e-4,1e-3);
        \coordinate (C) at (axis cs:6e-4,1e-2);

        \draw[thick] (A)--(B)--(C)--cycle;

        \node at (axis cs:4e-4,4e-4) {\tiny $1$};
        \node at (axis cs:1e-3,4e-3) {\tiny $2.0037$};

        \coordinate (D) at (axis cs:2e-3,1.5e2);
        \coordinate (E) at (axis cs:6e-3,1.5e2);
        \coordinate (F) at (axis cs:6e-3,5.0e2);

        \draw[thick] (D)--(E)--(F)--cycle;

        \node at (axis cs:4e-3,4e1) {\tiny $1$};
        \node at (axis cs:9.8e-3,3e2) {\tiny $1.0817$};

    \end{loglogaxis}
\end{tikzpicture}%
    }
    \caption{Convergence of state variables for cylindrical pair (PH-MP).}
    \label{fig:convergence_state}
\end{figure}

For the Lagrange multiplier $\vec{\lambda}$,
convergence reduces to first-order. This is consistent with the expected behavior of the method. For convergence of the state variables, both the PH-MP and the PH-MP-GGL integrator yield exactly the same results cf. Figure \ref{fig:convergence_state}. 

\subsection{Closed loop multibody system}\label{sec:closed_loop_MBS}
\begin{figure}[H]
    \centering
    \begin{subfigure}[t]{0.53\textwidth}
        \centering
        \def\svgwidth{\textwidth}%
\begingroup%
  \makeatletter%
  \providecommand\color[2][]{%
    \errmessage{(Inkscape) Color is used for the text in Inkscape, but the package 'color.sty' is not loaded}%
    \renewcommand\color[2][]{}%
  }%
  \providecommand\transparent[1]{%
    \errmessage{(Inkscape) Transparency is used (non-zero) for the text in Inkscape, but the package 'transparent.sty' is not loaded}%
    \renewcommand\transparent[1]{}%
  }%
  \providecommand\rotatebox[2]{#2}%
  \newcommand*\fsize{\dimexpr\f@size pt\relax}%
  \newcommand*\lineheight[1]{\fontsize{\fsize}{#1\fsize}\selectfont}%
  \ifx\svgwidth\undefined%
    \setlength{\unitlength}{299.52184242bp}%
    \ifx\svgscale\undefined%
      \relax%
    \else%
      \setlength{\unitlength}{\unitlength * \real{\svgscale}}%
    \fi%
  \else%
    \setlength{\unitlength}{\svgwidth}%
  \fi%
  \global\let\svgwidth\undefined%
  \global\let\svgscale\undefined%
  \makeatother%
  \begin{picture}(1,0.55940102)%
    \lineheight{1}%
    \setlength\tabcolsep{0pt}%
    \put(0,0){\includegraphics[width=\unitlength,page=1]{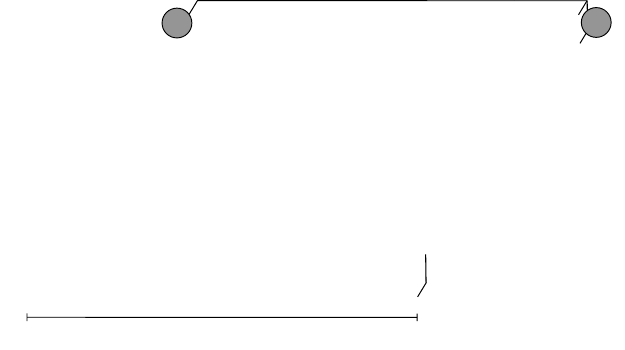}}%
    \put(0.35169945,0.00000001){\color[rgb]{0,0,0}\makebox(0,0)[lt]{\lineheight{1.25}\smash{\begin{tabular}[t]{l}$l$\end{tabular}}}}%
    \put(0.88050102,0.26517176){\color[rgb]{0,0,0}\makebox(0,0)[lt]{\lineheight{1.25}\smash{\begin{tabular}[t]{l}$\bm{F}$\end{tabular}}}}%
    \put(0.96528625,0.26704978){\color[rgb]{0,0,0}\makebox(0,0)[lt]{\lineheight{1.25}\smash{\begin{tabular}[t]{l}$\bm{\tau}$\end{tabular}}}}%
    \put(0.59050285,0.28404065){\color[rgb]{0,0,0}\makebox(0,0)[lt]{\lineheight{1.25}\smash{\begin{tabular}[t]{l}$\bm{x}$\end{tabular}}}}%
    \put(0.5630329,0.38950354){\color[rgb]{0,0,0}\makebox(0,0)[lt]{\lineheight{1.25}\smash{\begin{tabular}[t]{l}$\bm{y}$\end{tabular}}}}%
    \put(0.50830945,0.4200462){\color[rgb]{0,0,0}\makebox(0,0)[lt]{\lineheight{1.25}\smash{\begin{tabular}[t]{l}$\bm{z}$\end{tabular}}}}%
    \put(0,0){\includegraphics[width=\unitlength,page=2]{drawings/closed-loop-MBS_initial.pdf}}%
  \end{picture}%
\endgroup%
        \caption{Initial configuration of the closed loop multibody system.}
        \label{fig:MBS_sketch}
    \end{subfigure}
    \hfill
    \begin{subfigure}[t]{0.38\textwidth}
        \centering
        \def\svgwidth{\textwidth}%
\begingroup%
  \makeatletter%
  \providecommand\color[2][]{%
    \errmessage{(Inkscape) Color is used for the text in Inkscape, but the package 'color.sty' is not loaded}%
    \renewcommand\color[2][]{}%
  }%
  \providecommand\transparent[1]{%
    \errmessage{(Inkscape) Transparency is used (non-zero) for the text in Inkscape, but the package 'transparent.sty' is not loaded}%
    \renewcommand\transparent[1]{}%
  }%
  \providecommand\rotatebox[2]{#2}%
  \newcommand*\fsize{\dimexpr\f@size pt\relax}%
  \newcommand*\lineheight[1]{\fontsize{\fsize}{#1\fsize}\selectfont}%
  \ifx\svgwidth\undefined%
    \setlength{\unitlength}{178.15680556bp}%
    \ifx\svgscale\undefined%
      \relax%
    \else%
      \setlength{\unitlength}{\unitlength * \real{\svgscale}}%
    \fi%
  \else%
    \setlength{\unitlength}{\svgwidth}%
  \fi%
  \global\let\svgwidth\undefined%
  \global\let\svgscale\undefined%
  \makeatother%
  \begin{picture}(1,0.81204603)%
    \lineheight{1}%
    \setlength\tabcolsep{0pt}%
    \put(0,0){\includegraphics[width=\unitlength,page=1]{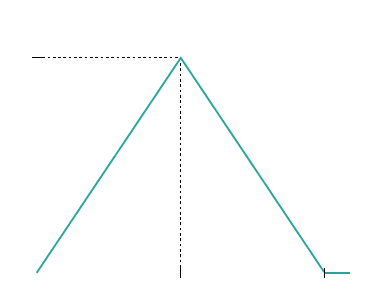}}%
    \put(-0.00167733,0.63608188){\color[rgb]{0,0,0}\makebox(0,0)[lt]{\lineheight{1.25}\smash{\begin{tabular}[t]{l}$f_m$\end{tabular}}}}%
    \put(0.11696657,0.77341806){\color[rgb]{0,0,0}\makebox(0,0)[lt]{\lineheight{1.25}\smash{\begin{tabular}[t]{l}$f(t)$\end{tabular}}}}%
    \put(0.98113,0.01163354){\color[rgb]{0,0,0}\makebox(0,0)[lt]{\lineheight{1.25}\smash{\begin{tabular}[t]{l}$t$\end{tabular}}}}%
    \put(0.85541337,0.00069063){\color[rgb]{0,0,0}\makebox(0,0)[lt]{\lineheight{1.25}\smash{\begin{tabular}[t]{l}1\end{tabular}}}}%
    \put(0.44232907,0.00076468){\color[rgb]{0,0,0}\makebox(0,0)[lt]{\lineheight{1.25}\smash{\begin{tabular}[t]{l}0.5\end{tabular}}}}%
    \put(0,0){\includegraphics[width=\unitlength,page=2]{drawings/external_load_MBS.pdf}}%
  \end{picture}%
\endgroup%
        \caption{Function of external loads.}
        \label{fig:MBS_loads}
    \end{subfigure}
    \caption{Closed loop multibody system.}\label{fig:closed_loop_MBS_sketch}
\end{figure}

This benchmark problem, adapted from \cite{betsch2002,kinon2024}, is used to assess the ability of the proposed redundant-coordinate formulation to handle closed kinematic loops under transient external loads, while preserving
physical invariants after load removal, as discussed in Section \ref{sec:numerical_integration}.
\begin{table}[H]
    \centering
    \caption{Simulation parameters for the closed loop multibody system.}
    \label{tab:closed_loop_MBS_setup}
    \renewcommand{\arraystretch}{1.3}
    \begin{tabular}{cccccccccc}
        \hline
        $h$ & $t_\text{end}$ & $l$ & $A$ & $\rho$ & $m$        & $f_m$ \\
        \hline
        0.1 & 10             & 10  & 1   & 1      & $\rho A l$ & 100   \\
        \hline
    \end{tabular}
\end{table}
The system consists of four rigid bars of length $l$ and square cross-section $A$, interconnected by four spherical joints to form a closed loop. Simulation parameters are given in Table \ref{tab:closed_loop_MBS_setup}.
Figure \ref{fig:MBS_sketch} shows a sketch of the initial configuration of the system.

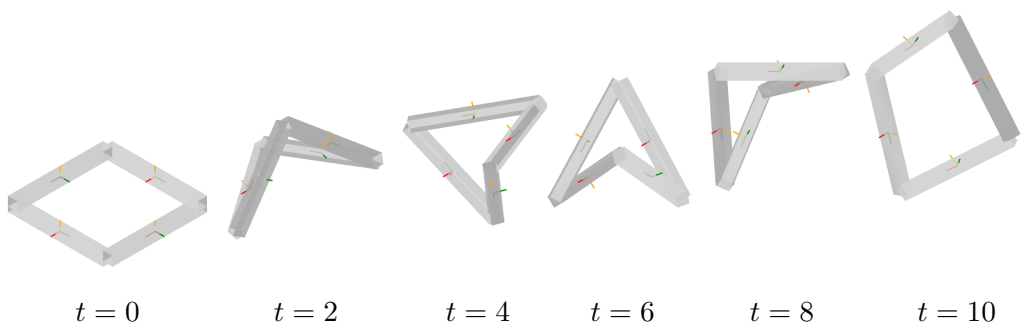
\begin{figure}[H]
    \centering
    \def\svgwidth{0.85\textwidth}%
\begingroup%
  \makeatletter%
  \providecommand\color[2][]{%
    \errmessage{(Inkscape) Color is used for the text in Inkscape, but the package 'color.sty' is not loaded}%
    \renewcommand\color[2][]{}%
  }%
  \providecommand\transparent[1]{%
    \errmessage{(Inkscape) Transparency is used (non-zero) for the text in Inkscape, but the package 'transparent.sty' is not loaded}%
    \renewcommand\transparent[1]{}%
  }%
  \providecommand\rotatebox[2]{#2}%
  \newcommand*\fsize{\dimexpr\f@size pt\relax}%
  \newcommand*\lineheight[1]{\fontsize{\fsize}{#1\fsize}\selectfont}%
  \ifx\svgwidth\undefined%
    \setlength{\unitlength}{1208.62767017bp}%
    \ifx\svgscale\undefined%
      \relax%
    \else%
      \setlength{\unitlength}{\unitlength * \real{\svgscale}}%
    \fi%
  \else%
    \setlength{\unitlength}{\svgwidth}%
  \fi%
  \global\let\svgwidth\undefined%
  \global\let\svgscale\undefined%
  \makeatother%
  \begin{picture}(1,0.30999331)%
    \lineheight{1}%
    \setlength\tabcolsep{0pt}%
    \put(0,0){\includegraphics[width=\unitlength,page=1]{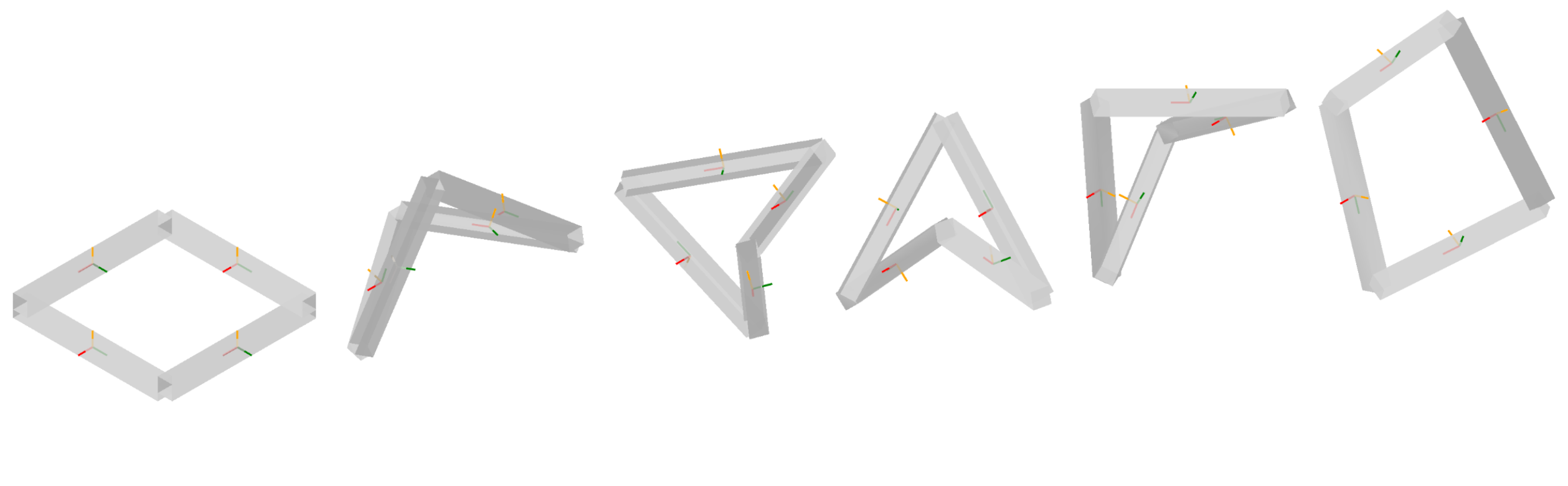}}%
    \put(0.07367499,0.000){\color[rgb]{0,0,0}\makebox(0,0)[lt]{\lineheight{1.25}\smash{\begin{tabular}[t]{l}$t=0$\end{tabular}}}}%
    \put(0.26662141,0.000){\color[rgb]{0,0,0}\makebox(0,0)[lt]{\lineheight{1.25}\smash{\begin{tabular}[t]{l}$t=2$\end{tabular}}}}%
    \put(0.43377973,0.000){\color[rgb]{0,0,0}\makebox(0,0)[lt]{\lineheight{1.25}\smash{\begin{tabular}[t]{l}$t=4$\end{tabular}}}}%
    \put(0.57446266,0.000){\color[rgb]{0,0,0}\makebox(0,0)[lt]{\lineheight{1.25}\smash{\begin{tabular}[t]{l}$t=6$\end{tabular}}}}%
    \put(0.72952742,0.000){\color[rgb]{0,0,0}\makebox(0,0)[lt]{\lineheight{1.25}\smash{\begin{tabular}[t]{l}$t=8$\end{tabular}}}}%
    \put(0.89239045,0.000){\color[rgb]{0,0,0}\makebox(0,0)[lt]{\lineheight{1.25}\smash{\begin{tabular}[t]{l}$t=10$\end{tabular}}}}%
  \end{picture}%
\endgroup%
    \vspace{0.5cm}
    \caption{Motion of system over simulation time. }
    \label{fig:MBS_motion}
    \vspace{-0.5cm}
\end{figure}
\begin{figure}[H]
    \centering
    \begin{subfigure}[b]{0.48\textwidth}
        \centering
        \adjustbox{valign=b}{%
            \definecolor{mycolor1}{rgb}{0.90196,0.62353,0.00000}%
\definecolor{mycolor2}{rgb}{0.33725,0.70588,0.91373}%
\definecolor{mycolor3}{rgb}{0.80000,0.47451,0.65490}%

\begin{tikzpicture}
    \begin{axis}[%
            width=0.65\textwidth,
            height=0.1\textheight,
            scale only axis,
            xmin=0,
            xmax=10,
            ymin=-10,
            ymax=2200,
            xtick distance = 1,
            unbounded coords=discard,
            xlabel={$t^n$},
            ylabel={$H(\vec{x}\n)$},
            xmajorgrids,
            ymajorgrids,
            axis background/.style={fill=white},
            legend style={at={(0.97,0.03)}, anchor=south east, draw=white!15!black, fill=white, font=\footnotesize},
            legend entries={PH-MP, Livens-EM}
        ]

        \addplot [color=kit-royalblue100, line width=1.5pt, solid, smooth]
        table[col sep=comma, x=time, y=system_total_energy_current_time] {data/e108_closed_loop_MBS_directors_ph_results.csv};

        \addplot [color=kit-green100, dashed, line width=1pt]
        table[row sep=crcr]{%
                0	0\\
                0.1	0.880082\\
                0.2	14.0589\\
                0.3	70.9105\\
                0.4	222.664\\
                0.5	538.597\\
                0.6	986.583\\
                0.7	1423.55\\
                0.8	1781.86\\
                0.9	2014.93\\
                1	2095.48\\
                1.1	2095.48\\
                1.2	2095.48\\
                1.3	2095.48\\
                1.4	2095.48\\
                1.5	2095.48\\
                1.6	2095.48\\
                1.7	2095.48\\
                1.8	2095.48\\
                1.9	2095.48\\
                2	2095.48\\
                2.1	2095.48\\
                2.2	2095.48\\
                2.3	2095.48\\
                2.4	2095.48\\
                2.5	2095.48\\
                2.6	2095.48\\
                2.7	2095.48\\
                2.8	2095.48\\
                2.9	2095.48\\
                3	2095.48\\
                3.1	2095.48\\
                3.2	2095.48\\
                3.3	2095.48\\
                3.4	2095.48\\
                3.5	2095.48\\
                3.6	2095.48\\
                3.7	2095.48\\
                3.8	2095.48\\
                3.9	2095.48\\
                4	2095.48\\
                4.1	2095.48\\
                4.2	2095.48\\
                4.3	2095.48\\
                4.4	2095.48\\
                4.5	2095.48\\
                4.6	2095.48\\
                4.7	2095.48\\
                4.8	2095.48\\
                4.9	2095.48\\
                5	2095.48\\
                5.1	2095.48\\
                5.2	2095.48\\
                5.3	2095.48\\
                5.4	2095.48\\
                5.5	2095.48\\
                5.6	2095.48\\
                5.7	2095.48\\
                5.8	2095.48\\
                5.9	2095.48\\
                6	2095.48\\
                6.1	2095.48\\
                6.2	2095.48\\
                6.3	2095.48\\
                6.4	2095.48\\
                6.5	2095.48\\
                6.6	2095.48\\
                6.7	2095.48\\
                6.8	2095.48\\
                6.9	2095.48\\
                7	2095.48\\
                7.1	2095.48\\
                7.2	2095.48\\
                7.3	2095.48\\
                7.4	2095.48\\
                7.5	2095.48\\
                7.6	2095.48\\
                7.7	2095.48\\
                7.8	2095.48\\
                7.9	2095.48\\
                8	2095.48\\
                8.1	2095.48\\
                8.2	2095.48\\
                8.3	2095.48\\
                8.4	2095.48\\
                8.5	2095.48\\
                8.6	2095.48\\
                8.7	2095.48\\
                8.8	2095.48\\
                8.9	2095.48\\
                9	2095.48\\
                9.1	2095.48\\
                9.2	2095.48\\
                9.3	2095.48\\
                9.4	2095.48\\
                9.5	2095.48\\
                9.6	2095.48\\
                9.7	2095.48\\
                9.8	2095.48\\
                9.9	2095.48\\
                10	2095.48\\
            };

    \end{axis}
\end{tikzpicture}
        }
    \end{subfigure}
    \hfill
    \begin{subfigure}[b]{0.48\textwidth}
        \centering
        \adjustbox{valign=b}{%
            \definecolor{mycolor1}{rgb}{0.90196,0.62353,0.00000}%
\definecolor{mycolor2}{rgb}{0.33725,0.70588,0.91373}%
\definecolor{mycolor3}{rgb}{0.80000,0.47451,0.65490}%

\begin{tikzpicture}
    \begin{axis}[%
            width=0.65\textwidth,
            height=0.1\textheight,
            scale only axis,
            xmin=1,
            xmax=10,
            ymin = -2e-12,
            ymax = 2e-12,
            ytick ={ -2e-12, -1e-12, 0, 1e-12, 2e-12},
            xtick ={1,10},
            xlabel={$t^n$},
            ylabel={$H(\bm{x}^{n+1})-H(\bm{x}^{n})$},
            xmajorgrids,
            ymajorgrids,
            axis background/.style={fill=white},
        ]

        \addplot [color=kit-royalblue100, line width=1.0pt, solid, smooth]
        table[col sep=comma, x=time, y=system_total_energy_interval_increment] {data/e108_closed_loop_MBS_directors_ph_results_energy_increment.csv};

    \end{axis}
\end{tikzpicture}
        }
    \end{subfigure}
    \caption{Hamiltonian for closed loop multibody system (using PH-MP).}\label{fig:energy_closed_loop_MBS}
\end{figure}
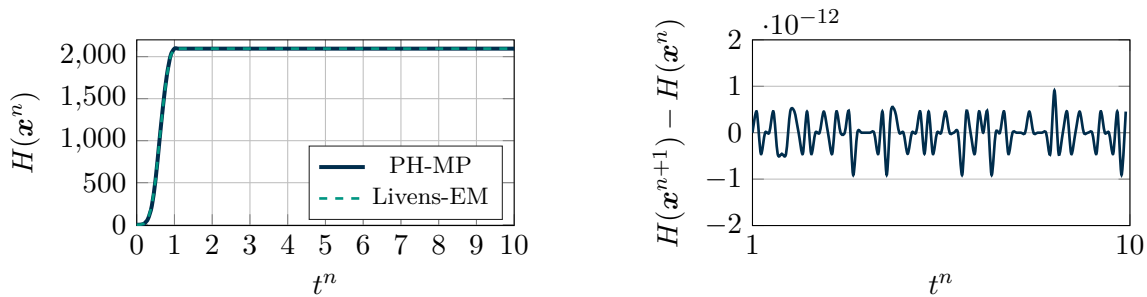
\begin{figure}[H]
    \begin{subfigure}[b]{0.48\textwidth}
        \centering
        \adjustbox{valign=b}{%
            \definecolor{mycolor1}{rgb}{0.90196,0.62353,0.00000}%
\definecolor{mycolor2}{rgb}{0.33725,0.70588,0.91373}%
\definecolor{mycolor3}{rgb}{0.80000,0.47451,0.65490}%

\begin{tikzpicture}
    \begin{axis}[%
            width=0.65\textwidth,
            height=0.10\textheight,
            scaled y ticks = true,
            scale only axis,
            xmin=0,
            xmax=10,
            xtick distance = 1,
            ymin=-20,
            ymax=320,
            unbounded coords=discard,
            xlabel={$t^n$},
            ylabel={$\bm{L}(\vec{x}\n)$},
            xmajorgrids,
            ymajorgrids,
            axis background/.style={fill=white},
            legend style={at={(0.97,0.03)}, anchor=south east, draw=white!15!black, fill=white, font=\footnotesize},
            legend cell align=left,
            legend entries={$L_x$, $L_y$, $L_z$}
        ]

        \addplot [color=kit-green100, line width=1.5pt, solid, smooth]
        table[col sep=comma, x=time, y=system_angular_momentum_x_current_time] {data/e108_closed_loop_MBS_directors_ph_results.csv};

        \addplot [color=kit-royalblue100, line width=1.5pt, dotted, smooth]
        table[col sep=comma, x=time, y=system_angular_momentum_y_current_time] {data/e108_closed_loop_MBS_directors_ph_results.csv};

        \addplot [color=kit-royalblue50, line width=1.5pt, dashed, smooth]
        table[col sep=comma, x=time, y=system_angular_momentum_z_current_time] {data/e108_closed_loop_MBS_directors_ph_results.csv};

    \end{axis}
\end{tikzpicture}
        }
    \end{subfigure}
    \hfill
    \begin{subfigure}[b]{0.48\textwidth}
        \centering
        \adjustbox{valign=b}{%
            \definecolor{mycolor1}{rgb}{0.90196,0.62353,0.00000}%
\definecolor{mycolor2}{rgb}{0.33725,0.70588,0.91373}%
\definecolor{mycolor3}{rgb}{0.80000,0.47451,0.65490}%

\begin{tikzpicture}
    \begin{axis}[%
            width=0.65\textwidth,
            height=0.10\textheight,
            scale only axis,
            xmin=1,
            xmax=10,
            xtick ={1,10},
            unbounded coords=discard,
            xlabel={$t^n$},
            ylabel={$\bm{L}(\bm{x}^{n+1})-\bm{L}(\bm{x}^{n})$},
            xmajorgrids,
            ymajorgrids,
            axis background/.style={fill=white},
            legend style={at={(0.97,0.03)}, anchor=south east, draw=white!15!black, fill=white, font=\footnotesize},
            legend cell align=left,
            legend entries={$L_x$, $L_y$, $L_z$}
        ]

        \addplot [color=kit-green100, line width=1.0pt, solid, smooth]
        table[col sep=comma, x=time, y=system_increment_angular_momentum_x] {data/e108_closed_loop_MBS_directors_ph_results_energy_increment.csv};

        \addplot [color=kit-royalblue100, line width=1.0pt, solid, smooth]
        table[col sep=comma, x=time, y=system_increment_angular_momentum_y] {data/e108_closed_loop_MBS_directors_ph_results_energy_increment.csv};

        \addplot [color=kit-royalblue50, line width=1.0pt, solid, smooth]
        table[col sep=comma, x=time, y=system_increment_angular_momentum_z] {data/e108_closed_loop_MBS_directors_ph_results_energy_increment.csv};

    \end{axis}
\end{tikzpicture}
        }
    \end{subfigure}
    \caption{Total angular momentum for closed loop multibody system (using PH-MP).}
    \label{fig:angular_momentum_MBS}
\end{figure}
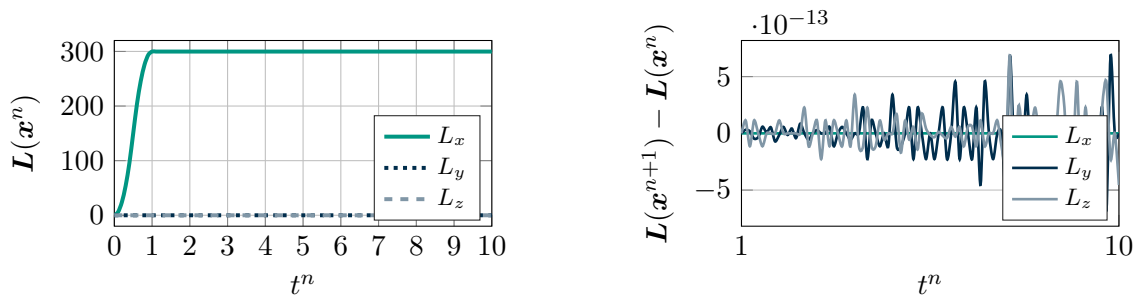

During $t \in [0,1]$, the center of mass of the first bar is subject to a ramp and decay load in the global $x$-direction (see Figure \ref{fig:MBS_loads}), which is given by
\begin{equation*}
    \vec{F}=8 f(t) \vec{e}_1, \quad \vec{\tau}=6 f(t) \vec{e}_1 \qquad \text{with} \qquad f(t) = \left\{\begin{aligned} 2 f_m t     & \quad \text{for}\: t \in  \left[0,0.5\right], \\
                2 f_m (1-t) & \quad \text{for}\: t \in  (0.5,1],            \\
                0           & \quad \text{for}\: t \in  (1,t_{\text{end}}].
    \end{aligned} \right.
\end{equation*}
The applied load induces both a global translation along the $x$-axis and a folding and unfolding motion of the structure, which persists after load removal. Figure \ref{fig:MBS_motion} presents
snapshots of this motion. As in the previous example, the total energy after load removal remains constant up to numerical round-off, see Figure \ref{fig:energy_closed_loop_MBS}. Further, it can be observed that the present method perfectly matches the result obtained with a unit-quaternion based energy-momentum consistent method \enquote{Livens-EM} from \cite{kinon_2024_conserving}.
This is also true for angular momentum, as shown in Figure \ref{fig:angular_momentum_MBS}. It should be noted that only the $L_x$ component remains nonzero, consistent with the loading direction.
This example confirms that the proposed formulation preserves invariant energy and angular momentum in closed-loop multibody dynamics under transient loading.
\subsection{Spatial slider-crank mechanism}\label{sec:slider_crank_mechanism}
\begin{figure}[H]
    \centering
    \def\svgwidth{0.65\textwidth}%
    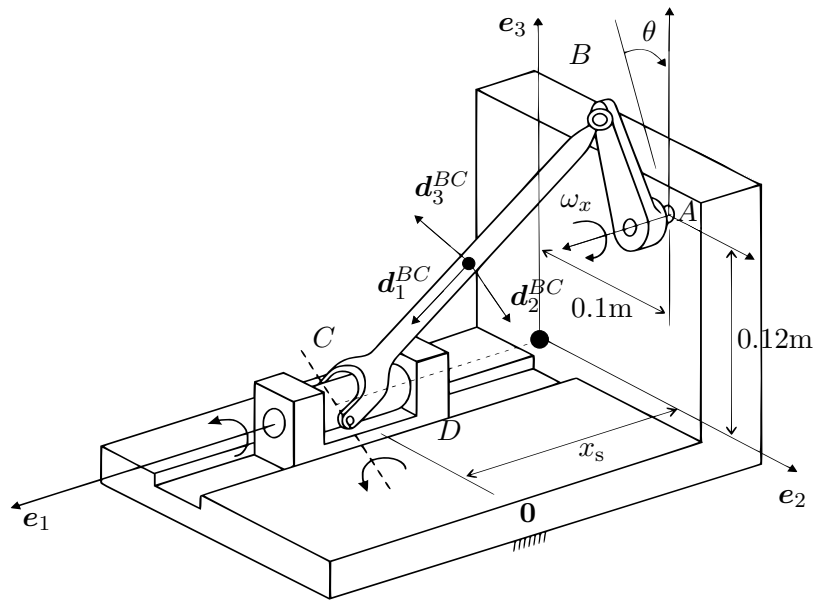%
    \vspace{0.5cm}
    \caption{Set-up spatial slider-crank mechanism. Figure adapted from \cite{masoudi2025}.}
    \label{fig:slider_crank}
    \vspace{-0.5cm}
\end{figure}

As a final benchmark, we apply the proposed framework to a spatial slider–crank mechanism, adapted from \cite{haug2021} and previously studied in \cite{masoudi2025, zhou2022}.
This example features a closed kinematic loop, fixed supports, and gravity $\vec{g}=-9.81 \vec{e}_3$.
\begin{table}[ht]
    \centering
    \caption{Simulation parameters of the spatial slider–crank benchmark.}
    \label{tab:spatial_slider_crank_setup}
    \renewcommand{\arraystretch}{1.3}
    \small
    \begin{tabular}{cccccccccccc}
        \hline
        $h$         & $t_\text{end}$ & $\omega_{x}(t=0)$ & $\theta_0$ & $l_{AB}$ & $a_{AB}$ & $m_{AB}$ & $l_{BC}$ & $a_{BC}$ & $m_{BC}$ & $a^{\text{block}}$ & $m^{\text{block}}$ \\
        \hline
        0.01, 0.001 & 5              & 6                 & 0          & 0.08     & 0.01     & 0.12     & 0.3      & 0.01     & 0.5      & 0.03               & 2.0                \\
        \hline
    \end{tabular}
\end{table}
The system is depicted in Figure~\ref{fig:slider_crank}  and consists of a crank $AB$ with mass $m_{AB}$, length $l_{AB}$ and a square cross-section characterized by side length $a_{AB}$, a connecting rod $BC$ with mass $m_{BC}$, length $l_{BC}$, and its square cross section characterized by $a_{BC}$, and a sliding block with mass $m^{\text{block}}$ modeled as a cube with side length $a^{\text{block}}$. The crank is connected to the wall by a revolute joint at $A$ and rotates with an initial
angular velocity $\omega_{x,0}$ at $t=0$ about the global $x$-axis. The block is constrained to the ground by a prismatic joint at $D$, allowing displacement along the $x$-axis. A spherical joint at $B$ and a universal joint
at $C$ connect the link to the crank and slider, respectively.

No further load is applied to the system. Key simulation parameters are summarized in
Table~\ref{tab:spatial_slider_crank_setup}.

In the simulated motion, the crank $AB$ rotates continuously about the global $\vec{e}_1$-axis. This rotation drives the connecting rod $BC$, which undergoes both translation and rotation.
As a consequence, the slider block moves along the $\vec{e}_1$-axis.

\begin{figure}[H]
    \centering
    \definecolor{mycolor1}{rgb}{0.90196,0.62353,0.00000}%
\definecolor{mycolor2}{rgb}{0.33725,0.70588,0.91373}%
\definecolor{mycolor3}{rgb}{0.80000,0.47451,0.65490}%

\begin{tikzpicture}
    \begin{axis}[%
            width=0.84\textwidth,
            height=0.15\textheight,
            scale only axis,
            xmin=0,
            xmax=5.0,
            unbounded coords=discard,
            xlabel={$t^n$},
            ylabel={$x_\mathrm{s}\n$},
            xmajorgrids,
            ymajorgrids,
            axis background/.style={fill=white},
            legend columns = -1,
            legend style={at={(1,1.03)}, anchor=south east, draw=white!15!black, 
            fill=white, font=\footnotesize},
            legend cell align=left,
            legend entries={PH-MP for $h=0.01$, Rosenbrock for $h=\frac{1}{60}$ \cite{masoudi2025}, Lie group generalized-$\alpha$ \cite{zhou2022}}
        ]

        \addplot [color=kit-green100, line width=1.5pt, solid, smooth]
        table[col sep=semicolon, x index=0, y index=1] {data/slider_crank_results.csv};

        \addplot [color=kit-royalblue100, line width=1.5pt, densely dashed, smooth]
        table[col sep=semicolon, x index=7, y index=8] {data/slider_crank_results.csv};

        \addplot [color=kit-royalblue50, line width=1.5pt, densely dotted, smooth]
        table[col sep=semicolon,x index=10, y index=11] {data/slider_crank_results.csv};

    \end{axis}
\end{tikzpicture}
    \caption{Comparison of simulated slider displacement with reference results.}
    \label{fig:x-pos-slider-crank}
    \vspace{-0.5cm}
\end{figure}
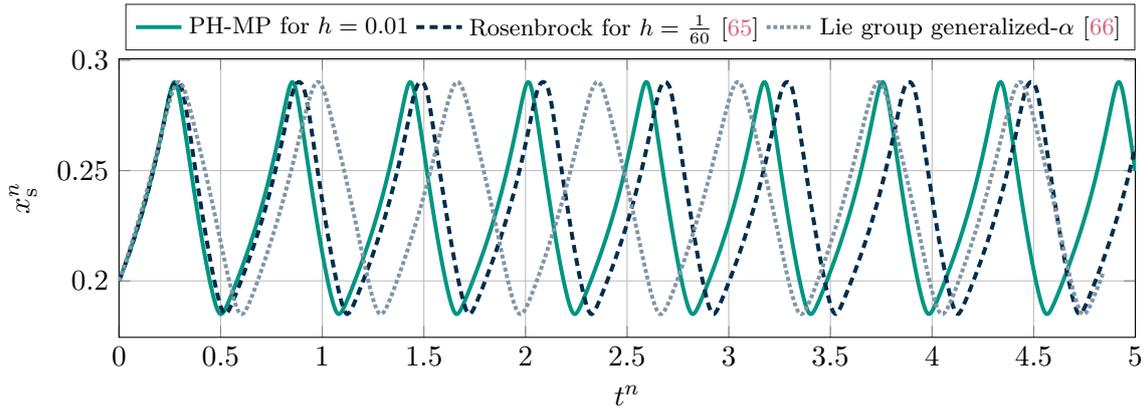

Figure~\ref{fig:x-pos-slider-crank} compares the resulting slider displacement $x_\mathrm{s}$ computed with our method (\enquote{PH-MP}) with reference results from \cite{masoudi2025, zhou2022}, which have been obtained from the database \cite{masoudi2025}. While the amplitudes are in good agreement, phase shifts are visible. 
Reducing the time step by two orders of magnitude to 
$h=0.0001$ did not produce any significant changes in the position profile obtained with our integrator, indicating that the result for $h=0.01$ is already close to the converged solution.
Both the Rosenbrock integrator in \cite{masoudi2025} and the Lie group generalized-$\alpha$ method in \cite{zhou2022} show a slower response. Since all system parameters are identical,
these differences likely stem from the different integration schemes.


Regarding energy conservation, the authors of \cite{masoudi2025} do not report energy results, while \cite{zhou2022} observes a gradual numerical damping, which is conforming with the observed phase shifts.

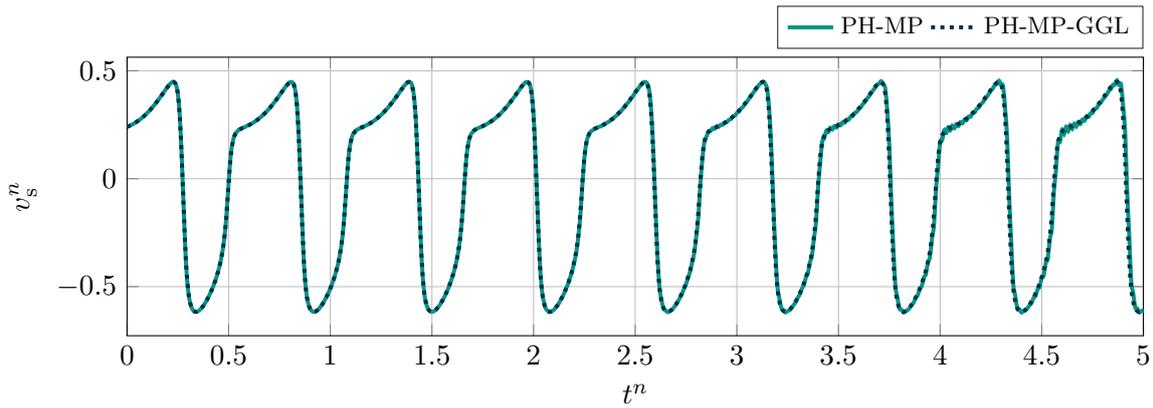
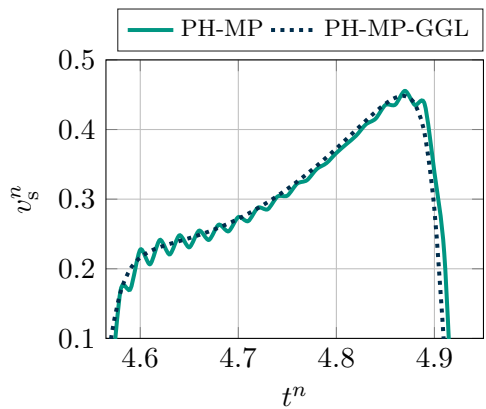
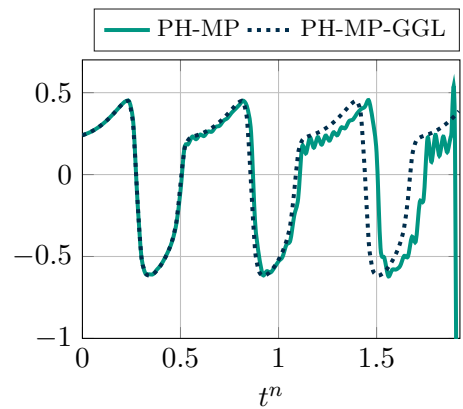
\begin{figure}[H]
    \centering
    \begin{subfigure}[t]{1\textwidth}
        \centering
        \adjustbox{valign=t}{%
            \definecolor{mycolor1}{rgb}{0.90196,0.62353,0.00000}%
\definecolor{mycolor2}{rgb}{0.33725,0.70588,0.91373}%
\definecolor{mycolor3}{rgb}{0.80000,0.47451,0.65490}%

\begin{tikzpicture}
    \begin{axis}[%
            width=0.84\textwidth,
            height=0.15\textheight,
            scale only axis,
            xmin=0,
            xmax=5.0,
            unbounded coords=discard,
            xlabel={$t^n$},
            ylabel={$v_\mathrm{s}\n$},
            xmajorgrids,
            ymajorgrids,
            axis background/.style={fill=white},
            legend columns = -1,
            legend style={at={(1,1.03)}, anchor=south east, draw=white!15!black, fill=white, font=\footnotesize},
            legend cell align=left,
            legend entries={PH-MP, PH-MP-GGL}
        ]

        \addplot [color=kit-green100, line width=1.5pt, solid, smooth]
        table[col sep=semicolon, x index=0, y index=3] {data/slider_crank_results.csv};


        \addplot [color=kit-royalblue100, line width=1.5pt, dotted, smooth]
        table[col sep=semicolon, x index=0, y index=4] {data/slider_crank_results.csv};

    \end{axis}
\end{tikzpicture}}
        \caption{Velocity profile for $h=0.01$.}
        \label{fig:x-vel-slider-crank}
    \end{subfigure}
    \vfill
    \vspace{5mm}
    \begin{subfigure}[t]{0.48\textwidth}
        \centering
        \adjustbox{valign=t}{%
            \definecolor{mycolor1}{rgb}{0.90196,0.62353,0.00000}%
\definecolor{mycolor2}{rgb}{0.33725,0.70588,0.91373}%
\definecolor{mycolor3}{rgb}{0.80000,0.47451,0.65490}%

\begin{tikzpicture}
    \begin{axis}[%
            width=0.65\textwidth,
            height=0.15\textheight,
            scale only axis,
            xmin=4.565,
            xmax=4.95,
            ymin=0.1,
            ymax=0.5,
            unbounded coords=discard,
            xlabel={$t^n$},
            ylabel={$v_\mathrm{s}\n$},
            xmajorgrids,
            ymajorgrids,
            axis background/.style={fill=white},
            legend columns = -1,
            legend style={at={(1,1.03)}, anchor=south east, draw=white!15!black, fill=white, font=\footnotesize},
            legend cell align=left,
            legend entries={PH-MP, PH-MP-GGL}
        ]

        \addplot [color=kit-green100, line width=1.5pt, solid, smooth]
        table[col sep=semicolon, x index=0, y index=3] {data/slider_crank_results.csv};


        \addplot [color=kit-royalblue100, line width=1.5pt, dotted, smooth]
        table[col sep=semicolon, x index=0, y index=4] {data/slider_crank_results.csv};

    \end{axis}
\end{tikzpicture}
        }
        \caption{Zoomed-in view for $h=0.01$.}
        \label{fig:x-vel-slider-crank-zoom}
    \end{subfigure}
    \hfill
    \begin{subfigure}[t]{0.48\textwidth}
        \centering
        \adjustbox{valign=t}{%
            \definecolor{mycolor1}{rgb}{0.90196,0.62353,0.00000}%
\definecolor{mycolor2}{rgb}{0.33725,0.70588,0.91373}%
\definecolor{mycolor3}{rgb}{0.80000,0.47451,0.65490}%

\begin{tikzpicture}
    \begin{axis}[%
            width=0.65\textwidth,
            height=0.15\textheight,
            scale only axis,
            xmin=0,
            xmax=1.925,
            ymin=-1,
            ymax=0.7,
            unbounded coords=discard,
            xlabel={$t^n$},
            ylabel={$v_\mathrm{s}\n$},
            xmajorgrids,
            ymajorgrids,
            axis background/.style={fill=white},
            legend columns = -1,
            legend style={at={(1,1.03)}, anchor=south east, draw=white!15!black, fill=white, font=\footnotesize},
            legend cell align=left,
            legend entries={PH-MP, PH-MP-GGL}
        ]

        \addplot [color=kit-green100, line width=1.5pt, solid, smooth]
        table[col sep=semicolon, x index=2, y index=0] {data/slider_crank_results_h=0.02.csv};

        \addplot [color=kit-royalblue100, line width=1.5pt, dotted, smooth]
        table[col sep=semicolon, x index=2, y index=1] {data/slider_crank_results_h=0.02.csv};

    \end{axis}
\end{tikzpicture}
        }
        \caption{Velocity profile for $h=0.02$.}
        \label{fig:x-vel-slider-crank-zoom_h=0.02}
    \end{subfigure}
    \caption{Velocity profile for $x$-velocity of slider at different time step sizes.}
    \label{fig:x-vel-slider-crank_all}
\end{figure}

Finally, we analyze the velocity of the slider  $v_\mathrm{s}$ as shown in Figure \ref{fig:x-vel-slider-crank_all}.
For a time step size of $h=0.01$ using the PH-MP integrator (see the green curves), the velocity profile remains qualitatively correct but exhibits small oscillations toward the end of the simulation, see Figure \ref{fig:x-vel-slider-crank} for the complete simulation time, and Figure \ref{fig:x-vel-slider-crank-zoom} for a zoomed-in view.
These oscillations cannot be observed anymore when the step size is reduced to $h=0.001$. For a step size of $h=0.02$ however, the oscillations induced by the PH-MP method become so pronounced that the Newton Raphson method eventually fails to converge, as depicted in Figure \ref{fig:x-vel-slider-crank-zoom_h=0.02}.

Contrarily, by enforcing constraints not just on position, but also on velocity level, the PH-MP-GGL integrator (see blue dashed curve) does not exhibit any oscillations, see Figure \ref{fig:x-vel-slider-crank-zoom}, even in the case where the PH-MP integrator fails, see Figure \ref{fig:x-vel-slider-crank-zoom_h=0.02}.

Overall, the proposed PH-MP integrator reproduces the expected kinematics of the slider-crank mechanism nicely while ensuring strict energy conservation throughout the simulation. Extending the method to the PH-MP-GGL integrator offers yet another significant benefit of avoiding oscillations in the velocity profile. Crucially, this enables the use of substantially larger step sizes, thereby offsetting the additional computational cost associated with the exact enforcement of velocity-level constraints.

\section{Conclusion and outlook}\label{sec:conclusion}
This work has presented a systematic integration of the director formulation for rigid bodies into the PH framework, providing a robust and energy-consistent approach to modeling multibody systems.
By combining the singularity-free representation of large rotations from the director formulation with the modular structure of PH systems, we have successfully simulated complex multibody dynamics
while preserving fundamental conservation properties.
In addition to that, we included GGL-type stabilization in the present PH framework. In this way, we have successfully gained additional stability benefits and ensured constraint satisfaction on both position and velocity level without compromising the ability to conserve total energy and angular momentum.

These results directly address the contributions outlined in the introduction: \textbf{C1}, presenting a PH formulation for mechanical systems based on the Lagrange formulation; \textbf{C2}, specifying this formulation for multibody systems in the director framework, thereby ensuring a singularity-free, energy-consistent
description of large rotations with constant mass matrices; \textbf{C3}, ensuring consistency of energy, angular momentum and constraints in numerical simulations using structure-preserving discretization; \textbf{C4}, identifying the equivalence between kinematic pairs and PH interconnection matrices, thus bridging classical multibody theory and the PH interconnection approach; and \textbf{C5}, additionally satisfying velocity-level constraints and achieving index-reduction by using the GGL principle. Together, these contributions demonstrate the suitability of the proposed approach for simulating complex multibody systems in a modular and physically consistent way.

Looking ahead, several promising research directions emerge. Integrating control inputs via the PH ports could enable energy-shaping \cite{caasenbrood2022} and other advanced control strategies, exploiting the inherent passivity of the system.
Extending the framework to flexible multibody systems would allow for the simultaneous treatment of rigid and elastic components in a unified, energy-consistent manner. A first step in this direction has been made for planar geometrically exact beams in \cite{kinon2025}. Similarly, embedding multiphysics phenomena, such as
thermo-elasticity or electromechanical interactions, would capitalize on the modularity of the PH representation. Finally, the director formulation, while robust and instructive, introduces six additional constraints per rigid body, potentially increasing the computational cost. On the other hand, the director formulation comes with a sparse matrix formulation and a low degree of nonlinearity, including a constant mass matrix.
Alternatives such as quaternion-based representations could improve computational efficiency while maintaining the favorable conservation properties demonstrated in this work, however, at the cost of state-dependent mass matrices.

In summary, this work has demonstrated that the director formulation embedded within the PH framework provides a modular and physically consistent approach to modeling multibody systems.
By preserving energy and angular momentum, linking classical constraints with PH interconnections, and maintaining passivity in both continuous and discrete models, the proposed framework establishes a solid basis for simulating and controlling complex multibody systems and offers a foundation for future extensions to flexible, controlled or multiphysics applications.

\begin{acks}
PLK gratefully acknowledges funding by the Research Travel Grant of the Karlsruhe House of Young Scientists (KHYS). 
We are thankful for fruitful discussions with Simon R. Eugster (TU Eindhoven) and Manuel Schaller (TU Chemnitz).
\end{acks}
\begin{dci}
    The authors declare no known conflict of interest.
\end{dci}
\begin{code}
    The simulation results are publicly available under \cite{latussek2026}. 
The simulations have been conducted with the open-source dynamics toolkit \textit{pydykit} \cite{kinonbauerlatussek2026}.
\end{code}
\begin{authcontrib}
\vspace{-2em}
\begin{table*}[h]
    \centering
    \begin{tabular}{l p{14cm}}
        \textbf{LL}  &                    
        Methodology,
        software,
        validation,
        formal analysis,
        investigation,
        data curation,
        writing--original draft preparation,
        visualization.                    \\
        \textbf{PLK} & Conceptualization,
        methodology,
        software,
        validation,
        formal analysis,
        investigation,
        data curation,
        writing--original draft preparation,
        visualization,
        supervision.                      \\
        \textbf{PB}  &                    
        Methodology,
        validation,
        formal analysis,
        resources,
        writing--review and editing,
        supervision,
        project administration,
        funding acquisition.
    \end{tabular}
\end{table*}
\vspace{-1.5em}
\end{authcontrib}

\printbibliography

\appendix
\section*{Appendix}\label{annexsec:annex}

\section{Balance of angular momentum for the director rigid body} \label{app_ang_mom}

For ease of analysis, we consider an individual rigid body as described in \eqref{eq:ph_constraint_compact_simple} with total angular momentum \eqref{eq:angular_momentum_rigid_body}. Its time derivative reads
\begin{equation}
    \dot{\vec{L}} = \boldsymbol{\varphi} \times  m \ddot{\vec{\varphi}}  + \sum_{i=1}^{3} \vec{d}_i \times  E_i \ddot{\vec{d}}_i,
    \label{eq:balance_angular_momentum}
\end{equation}
due to the vanishing cross-product of a vector with itself. Considering the rigid body's equations of motion
\begin{subequations} \label{eq:multibody_system_equations_of_motion}
\begin{align}
        \dot{\boldsymbol{\varphi}}&= \bm{v}_{\boldsymbol{\varphi}}, \label{eq:evolution_position_trans_mbs}\\
        \dot{\bm{d}}_i  &= \bm{v}_{\bm{d}_i}, \label{eq:evolution_position_rot_mbs}\\
        m \dot{\bm{v}}_{\boldsymbol{\varphi}} &= -\nabla V(\boldsymbol{\varphi}) + \bm{f}_{\boldsymbol{\varphi}}, \label{eq:evolution_velocity_trans_mbs}\\
        E_i \dot{\bm{v}}_{\bm{d}_i} &= -\nabla\, \bm{g}_\text{d}(\bm{d}_i)^\top\,\bm{\lambda} \bm{f}_i, \label{eq:evolution_velocity_rot_mbs}\\
    \nabla \bm{g}(\bm{d}_i)\: \bm{v}_{\bm{d}_i} & = \bm{0}, \label{eq:director_constraints_mbs}
\end{align}
\end{subequations}
we insert \eqref{eq:evolution_velocity_trans_mbs} and \eqref{eq:evolution_velocity_rot_mbs} into \eqref{eq:balance_angular_momentum}. We further consider co-linearity of constraint forces and directors such that $\vec{f}^\mathrm{d}_i \times \vec{d}_i = \vec{0} $ and obtain
\begin{align}
    \dot{\vec{L}} =-\boldsymbol{\varphi} \times \nabla V(\boldsymbol{\varphi})+ \boldsymbol{\varphi}\times \vec{f}_{\boldsymbol{\varphi}}+ \vec{d}_i \times \vec{f}_i,
    \label{eq:balance_angular_momentum2}
\end{align}
showing that in the absence of external and potential forces, angular momentum is conserved. In the case of active gravitational forces, e.g. $\nabla V (\boldsymbol{\varphi})=- m g\, \vec{e}_3$, the component of the angular
momentum in the $\vec{e}_3$-direction remains constant, i.e. $\frac{\text{d}}{\text{d}t}(\vec{L}\transp \vec{e}_3)=0$.

In discrete time, the angular momentum balance can be shown to hold true for the discrete evolution equations \eqref{eq:port_hamiltonian_form_discrete} as well.
If applied to the rigid body formulation
the scheme is given by
\begin{subequations} \label{eq:discretized_eqs_full}
    \begin{align}
        \boldsymbol{\varphi}^{n+1} - \boldsymbol{\varphi}^n                                    & = h\, \vec{v}_{\boldsymbol{\varphi}}^{n+\frac{1}{2}}, \label{eq:discretized_eq_translational_com}                                                                                                             \\
        \vec{d}_i^{n+1} - \vec{d}_i^n                                                          & = h\, \vec{v}_{\vec{d}_i}^{n+\frac{1}{2}} \quad \text{for} \ i=1, 2, 3 , \label{eq:discretized_eq_rotational_com}                                                                                             \\
        m \left(\vec{v}_{\boldsymbol{\varphi}}^{n+1} - \vec{v}_{\boldsymbol{\varphi}}^n\right) & = -h\, \nabla V(\boldsymbol{\varphi}^{n+\frac{1}{2}})  +h\,\vec{f}_{\boldsymbol{\varphi}}^{n+\frac{1}{2}}, \label{eq:discretized_eq_translational_d}                                                          \\
        E_i \left(\vec{v}_{\vec{d}_i}^{n+1} - \vec{v}_{\vec{d}_i}^n\right)                     & = -h\, \nabla \vec{g}(\vec{d}_i^{n+\frac{1}{2}})\transp \vec{\lambda}_\text{d}^{n+\frac{1}{2}} + h\,\vec{f}_{\vec{d}_i}^{n+\frac{1}{2}} \quad \text{for} \ i=1, 2, 3 , \label{eq:discretized_eq_rotational_d} \\
        \nabla \vec{g} (\vec{d}_i^{n+\frac{1}{2}})\, \vec{v}_{\vec{d}_i}^{n+\frac{1}{2}}       & = \vec{0} .  \label{eq:discretized_eq_constraint}
    \end{align}
\end{subequations}

We now show that the discrete-time formulation under the midpoint rule also yields an angular momentum balance.
Evaluating the continuous angular momentum balance \eqref{eq:balance_angular_momentum}
at $t^n$ and $t^{n+1}$ yields
\begin{equation}
    \begin{aligned}
        \vec{L}(\vec{x}^{n+1})  - \vec{L}(\vec{x}^n) =  \vec{\varphi}^{n+1} \times m \vec{v}_{\boldsymbol{\varphi}}^{n+1} - \vec{\varphi}^{n} \times m\vec{v}_{\boldsymbol{\varphi}}^{n}  +  \vec{d}_i^{n+1} \times E_i \vec{v}_{\vec{d}_i}^{n+1} - \vec{d}_i^{n} \times E_i \vec{v}_{\vec{d}_i}^{n} .
    \end{aligned}
    \label{eq:change_angular_momentum}
\end{equation}
To evaluate each term, we use the standard identity
\begin{equation}
    \begin{aligned}
        \scriptscriptstyle \mathbf{a}^{n+1} \times \mathbf{b}^{n+1} - \mathbf{a}^n \times \mathbf{b}^n = \frac{1}{2}  \left[\left( \mathbf{a}^{n+1} - \mathbf{a}^n \right) \times  \left( \mathbf{b}^{n+1} + \mathbf{b}^n \right) +  \left( \mathbf{a}^{n+1} + \mathbf{a}^n \right) \times \left( \mathbf{b}^{n+1} - \mathbf{b}^n \right)\right] .
    \end{aligned}
    \label{eq:cross_product_difference}
\end{equation}
Applying \eqref{eq:cross_product_difference} to the translational term in \eqref{eq:change_angular_momentum} and substituting the discrete evolution laws \eqref{eq:discretized_eq_translational_com} and \eqref{eq:discretized_eq_translational_d} yields
\begin{equation}
    {\scriptscriptstyle
        \vec{\varphi}^{n+1} \times m\vec{v}_{\boldsymbol{\varphi}}^{n+1} - \vec{\varphi}^{n} \times m\vec{v}_{\boldsymbol{\varphi}}^{n}=
        \frac{h}{2} \left[ \left( \vec{v}_{\boldsymbol{\varphi}}^{n+1} + \vec{v}_{\boldsymbol{\varphi}}^{n} \right) \times \left( m \vec{v}_{\boldsymbol{\varphi}}^{n+1} + m \vec{v}_{\boldsymbol{\varphi}}^{n} \right)
            - \left( \vec{\varphi}^{n+1} + \vec{\varphi}^{n} \right) \times \nabla V\left( \vec{\varphi}^{n+\frac{1}{2}} \right) + \left( \vec{\varphi}^{n+1} + \vec{\varphi}^{n} \right) \times \vec{f}_{\boldsymbol{\varphi}}^{n+\frac{1}{2}}\right] .}
    \label{eq:change_translational_momentum_2}
\end{equation}
The first term vanishes due to the collinearity of the velocity vector with itself. A similar expansion for the director term results in
\begin{equation}
    {\scriptscriptstyle
        \vec{d}_i^{n+1} \times E_i\, \vec{v}_{\vec{d}_i}^{n+1} - \vec{d}_i^{n} \times E_i \,\vec{v}_{\vec{d}_i}^{n}=
        h\left[\frac{1}{2} \left( \vec{v}_{\vec{d}_i}^{n+1} + \vec{v}_{\vec{d}_i}^{n} \right) \times  \left(  E_i \vec{v}_{\vec{d}_i}^{n+1} +  E_i \vec{v}_{\vec{d}_i}^{n} \right)
            -  \vec{d}_i^{n+\frac{1}{2}}  \times \nabla \vec{g}( \vec{d}_i^{n+\frac{1}{2}})\transp\vec{\lambda}^{n+\frac{1}{2}} + \vec{d}_i^{n+\frac{1}{2}}  \times \vec{f}_{\vec{d}_i}^{n+\frac{1}{2}}\right].}
    \label{eq:change_rotational_momentum_2}
\end{equation}
Again, the first term vanishes due to collinearity, and the second vanishes due to collinearity of constraint forces and directors.
Combining both contributions yields
\begin{equation}
    \vec{L}(\vec{x}^{n+1})  - \vec{L}(\vec{x}^n) = - h \left[ \vec{\varphi}^{n+\frac{1}{2}}  \times \nabla V\left( \vec{\varphi}^{n+\frac{1}{2}} \right)- \vec{\varphi}^{n+\frac{1}{2}} \times \vec{f}_{\boldsymbol{\varphi}}^{n+\frac{1}{2}}-\,\vec{d}_i^{n+\frac{1}{2}}  \times \vec{f}_{\vec{d}_i}^{n+\frac{1}{2}}\right].
    \label{discrete_angmom_balance}
\end{equation}
Thus, if no external potential or external loads act on the system, the right-hand side vanishes, and angular momentum is exactly conserved under the midpoint rule.

\end{document}